\documentclass[hidelinks,onefignum,onetabnum]{siamart220329}

\usepackage{lineno}
\usepackage[utf8x]{inputenc}
\usepackage[T1]{fontenc}
\usepackage[figuresright]{rotating}
\usepackage{xcolor}
\usepackage{pdfcomment}
\usepackage[shortlabels]{enumitem}
\usepackage[english]{babel}
\usepackage{graphicx} 
\usepackage{amsfonts,amssymb,mathtools}
\usepackage{amstext}
\usepackage{upgreek}
\usepackage{stmaryrd}
\usepackage[tight,footnotesize]{subfigure}
\usepackage[makeroom]{cancel}
\usepackage[draft]{changes}
\setaddedmarkup{\textcolor{blue}{#1}}
\setdeletedmarkup{\textcolor{red}{\sout{#1}}}
\usepackage{import}
\usepackage{tikz}
\usepackage{circuitikz}
\usepackage{pgfplots}
\usepgfplotslibrary{groupplots}
\pgfplotsset{width=10cm,height=8cm,compat=1.16}
\usetikzlibrary{external}
\usepackage{pgfplotstable}
\usepackage[per-mode=symbol]{siunitx}
\usepackage{mathtools}
\usepackage{soul}
\usepackage[normalem]{ulem}
\usepackage{color}
\usepackage{blindtext}
\usepackage{tcolorbox}
\usepackage{accents}
\usepackage{bbm}
\usepackage{algorithmicx}
\usepackage{algorithm}
\usepackage{algpseudocode}
\usepackage{float}
\usepackage{collectbox}
\usepackage{caption}
\usepackage{siunitx}
\usepackage{setspace}
\graphicspath{{./images/}}

\algnewcommand\algorithmicinput{\textbf{INPUT:}}
\algnewcommand\INPUT{\item[\algorithmicinput]}
\algnewcommand\algorithmicoutput{\textbf{OUTPUT:}}
\algnewcommand\OUTPUT{\item[\algorithmicoutput]}

\numberwithin{equation}{section}

\numberwithin{thm}{section}

\numberwithin{prob}{section}

\renewcommand{\vec}[1]{\mathbf{\MakeLowercase{#1}}}
\graphicspath{{figures/}}

\newcommand{\bs}{\mathbf}

\newcommand{\bx}{\mbox{\boldmath$x$}}
\newcommand{\by}{\mbox{\boldmath$y$}}

\newcommand{\bn}{\mbox{\boldmath$n$}}
\newcommand{\bu}{\mbox{\boldmath$u$}}
\newcommand{\bv}{\mbox{\boldmath$v$}}

\newcommand{\ba}{\mbox{\boldmath$a$}}
\newcommand{\bb}{\mbox{\boldmath$b$}}
\newcommand{\bc}{\mbox{\boldmath$c$}}
\newcommand{\bd}{\mbox{\boldmath$d$}}
\newcommand{\be}{\mbox{\boldmath$e$}}

\newcommand{\bh}{\mbox{\boldmath$h$}}
\newcommand{\bj}{\mbox{\boldmath$j$}}

\newcommand{\br}{\mbox{\boldmath$r$}}

\newcommand{\bA}{\mbox{\boldmath$A$}}

\newcommand{\bM}{\mbox{\boldmath$M$}}

\newcommand{\bV}{\mbox{\boldmath$V$}}

\newcommand{\bAA}{\mbox{\boldmath$A$}}
\newcommand{\bBB}{\mbox{\boldmath$B$}}

\newcommand{\bEE}{\mbox{\boldmath$E$}}
\newcommand{\bFF}{\boldsymbol{F}}

\newcommand{\bHH}{\boldsymbol{H}}
\newcommand{\bJJ}{\mbox{\boldmath$J$}}
\newcommand{\bKK}{\boldsymbol{K}}

\newcommand{\bMM}{\boldsymbol{M}}

\newcommand{\bRR}{\boldsymbol{R}}
\newcommand{\bSS}{\boldsymbol{S}}

\newcommand{\bbeta}{\mbox{\boldmath$\eta$} }

\newcommand{\bmu}{\mbox{\boldmath$\mu$} }

\newcommand{\bsigma}{\mbox{\boldmath$\sigma$} }

\newcommand{\bzero}{\boldsymbol{0} }

\newcommand{\LLtwo}[2][]{\boldsymbol{L}^2#1(#2)}

\newcommand{\Hone}[2][]{H^1#1(#2)}

\newcommand{\Hcurl}[2][]{\boldsymbol{H}#1(\CurlSymb;#2)}
\newcommand{\Hdiv}[2][]{\boldsymbol{H}#1(\DivSymb;#2)}

\newcommand{\GradSymb}{\vec{grad}}
\newcommand{\CurlSymb}{\vec{curl}}
\newcommand{\DivSymb}{\mathrm{div}}

\newcommand{\Curl}[2][]{\mathrm{\CurlSymb}{#1}\,{#2}}
\newcommand{\Grad}[2][]{\mathrm{\GradSymb}{#1}\,{#2}}
\newcommand{\Div}[2][]{\mathrm{\DivSymb}{#1}\,{#2}}

\ifpdf
  \DeclareGraphicsExtensions{.eps,.pdf,.png,.jpg}
\else
  \DeclareGraphicsExtensions{.eps}
\fi

\newsiamremark{remark}{Remark}
\newsiamremark{hypothesis}{Hypothesis}
\crefname{hypothesis}{Hypothesis}{Hypotheses}
\newsiamthm{claim}{Claim}

\headers{FE-HMM Method for Confined Eddy Current Problems}{I. Niyonzima et al.}

\title{Magnetic Field Conforming Multiscale Formulations for Locally-Confined Nonlinear Eddy Current Problems Using the FE-HMM Method\thanks{Submitted to the editors on October 18, 2024.
}}

\author{
  Innocent Niyonzima\thanks{Univ. Grenoble Alpes, CNRS, Grenoble INP, G2ELab, F-38000 Grenoble, France, \email{innocent.niyonzima@univ-grenoble-alpes.fr}, \email{Gerard.Meunier@g2elab.grenoble-inp.fr}, \email{olivier.chadebec@g2elab.grenoble-inp.fr}, \email{nicolas.galopin@g2elab.grenoble-inp.fr}.}
  \and Gérard Meunier\footnotemark[2]
  \and Antoine Marteau\thanks{School of Mathematics, Monash University, Clayton, Victoria, 3800, Australia, \email{Antoine.Marteau@monash.edu}.}
  \and Ruth V. Sabariego\thanks{KU Leuven, Dept. Electrical Engineering (ESAT), Campus EnergyVille, B-3600 Genk, Belgium, \email{ruth.sabariego@kuleuven.be}.}
  \and Olivier Chadebec\footnotemark[2]
  \and Nicolas Galopin\footnotemark[2]
  \and Christophe Geuzaine\thanks{University of Li\`{e}ge, Dept. Electrical Engineering and Computer Science, Montefiore Institute B28, B-4000 Li\`{e}ge, Belgium, \email{cgeuzaine@uliege.be}.}
}

\ifpdf
\hypersetup{
  pdftitle={An Example Article},
  pdfauthor={D. Doe, P. T. Frank, and J. E. Smith}
}
\fi

\begin{document}

\maketitle

\begin{abstract}
    Magnetic composites used for the conversion of electrical energy often incorporate ferromagnetic inclusions insulated from each other to mitigate eddy current losses. Numerical models for these composites must be robust enough to address potential convergence issues arising from the presence of nonlinear magnetic inclusions and presence of significant confined eddy currents in the cell. This paper introduces an $\bh$-conforming multiscale formulation for magnetic composites in a periodic framework. The proposed method uses the Heterogeneous Multiscale Method (HMM) with two mesoscale problems: a magnetoquasistatic problem for upscaling the homogenized magnetic flux density $\bBB_M$ and a magnetostatic problem for upscaling the macroscale incremental reluctivity $(\partial \bBB_M/\partial \bHH_M)$. Additionally, the method uses relaxed Newton--Raphson schemes at both macro and mesoscale levels to mitigate the well-known NR convergence issue linked to nonlinear BH constitutive laws in $\bh$-conforming formulations. The accuracy and performance of the formulation are evaluated using 2D and 3D idealized periodic soft magnetic composites with linear and nonlinear BH curves. Furthermore, the paper demonstrates that the magnetoquasistatic mesoscale problem can be replaced by a magnetostatic problem, as long as the magnetic power within the cell substantially exceeds the eddy current losses.
\end{abstract}

\begin{keywords}
    Multiscale modeling, computational homogenization, finite element method, nonlinear magnetic problems, eddy current problems, $\bh$-conforming formulations. 
\end{keywords}

\begin{MSCcodes}
    35K55, 35M10, 65M60, 65Y0520, 788A48, 78M10, 78M35, 78M40
\end{MSCcodes}

\section{Introduction}
Asymptotic homogenization and multiscale methods have become essential tools for solving physical problems with composite materials. These methods have been used in many areas of engineering sciences including computational mechanics \cite{fish-homogenization-97, fish-homogenization-13, nguyen-homogenization-14}, flow in porous media \cite{auriault-homogenization-05, allaire-homogenization-07}, and for the homogenization of the electromagnetic wave problems with and without the scale separation assumption \cite{Wellander-homogenization-03, ouchetto-homogenization-06, ouchetto-homogenization-07, ciarlet-hmm-17, ohlberger-hmm-18, farhat-homogenization-08, bossavit-homogenization-09, cherednichenko-homogenization-19}, just to name a few. 
These methods have also been applied to homogenize magnetostatic and magnetoquasistatic (MQS) multiscale problems, such as the laminated cores and windings in electrical machines and transformers, soft ferrites or the idealized soft magnetic composites in high frequency transformers~\cite{bottauscio-chm-08, meunier-chm-10, niyonzima-hmm-12, appino-stat-homog-12, niyonzima-hmm-13, bottauscio-chm-13, niyonzima-chm-14, ruuskanen-chm-24, marteau-phdthesis-24}. They are indispensable numerical tools for modeling problems with 3D printed magnetic materials and for optimizing the microstructure of these 3D printed materials.

The main challenges posed by the homogenization of Maxwell's equations are: 
\begin{enumerate}
  \item The presence of nonlinear materials such as ferromagnetic or superconducting materials.
  \item The 3D geometries that can lead to electromagnetic fields with complex distributions.
  \item The presence of composite materials that leads to multiscale fields and the lack of scale separation which can occur when either the skin depth $\delta$ or the electromagnetic wavelength $\lambda$ become comparable to the size of inclusions in the microstructure.  
  \item The presence of confined eddy currents which may necessitate well-adapted upscaling techniques for the fields to be homogenized.  
  \item The consideration of stochastic distributions of inclusions in the microsctructure.   
\end{enumerate}
In this paper, we address the first four points using idealized soft magnetic composites in a periodic setting to validate our formulations for low frequency electromagnetic problems. 

The scientific literature on asymptotic homogenization for low frequency electromagnetic problems is extensive, especially for 2D multiscale problems. Bottauscio et al. proposed multiscale methods for 2D and 3D nonlinear problems using the \emph{asymptotic expansion} method, the \emph{multiscale finite element method} (MsFEM) and the \emph{variational multiscale methods} (VMS)~\cite{bottauscio-chm-08, bottauscio-chm-13}. The authors were able to solve static and dynamic multiscale problems in periodic composites with eddy currents considered in the conducting subdomain and displacement currents considered in the non-conducting subdomain of the periodic cell. While the proposed methods provided very accurate results, they failed to yield accurate results for dynamic problems with locally-confined eddy currents. In~\cite{niyonzima-hmm-12, niyonzima-hmm-13, niyonzima-chm-14}, the authors used the HMM method to solve magnetoquasistatic problems in 2D laminated materials and in idealized soft magnetic composites. This method was recently extended to 3D nonlinear problems \cite{marteau-phdthesis-24} and coupled with the multiharmonic approach to solve for the steady state regime for 3D nonlinear multiscale magnetoquasistatic problems \cite{ruuskanen-chm-24}. While this approach allows to provide simulations for composites in the steady state, it fails to provide simulations of the transient regime. In~\cite{meunier_chm-08, meunier-chm-10}, the authors proposed the multiscale $\bb$- and $\bh$-conforming formulations for 2D problems including problems with windings. This approach allowed the homogenization of multiscale magnetoquasistatic problems with upscaling techniques well-adapted for problems with confined eddy current problems. However, the study was limited to linear problems. The recent research work by Sch\"{o}binger and Hollaus allowed to simulate nonlinear multiscale problem with eddy currents for laminated structures using another version of the MsFEM method~\cite{hollaus-msfem-18, schobinger-msfem-19, hollaus-msfem-20}. The approach was extended by including model order reduction methods. However, the method is inherently limited to laminated composites. Ren et al. proposed the homogenization of nonlinear problems with eddy currents using the local orthogonal decomposition (LOD) method derived from the variational multiscale method (VMS) with 2D applications~\cite{ren-lod-21}. The work was recently extended to topology and shape optimization of the inclusions to get tailored macroscale properties~\cite{ren-optim-21, ren-optim-22}. 

Homogenization methods other than asymptotic homogenization have also been used for low frequency electromagnetic problems. The mean-field approach was used for 3D problems considering nonlinear hysteretic laws~\cite{daniel-coupling-14, corcolle-mfh-21}, and later extended to eddy current problems~\cite{preault-mfh-14}. However, this method cannot correctly handle 3D non-linear problems with strongly confined eddy currents. A classical homogenization method has been developed for non-linear multiscale problems with eddy currents in the context of laminated cores and windings~\cite{gyselinck-homogenization-04, gyselinck-homogenization-06}. Later, the approach was extended to account for global eddy currents circulating in the stack of laminations due, for instance, to contact manufacturing defects in the magnetic stack of laminations \cite{gyselinck-homogenization-16}. However, the method is not well-adapted for composites with a periodic or non-periodic arbitrary cell.

In this paper, we propose an $\bh$-conforming multiscale formulation and its numerical implementation for solving 3D multiscale electromagnetic problems using the heterogeneous multiscale method (HMM)~\cite{e-hmm-03, abdulle-hmm-12,abdulle-hmm-16}. In the context of the eddy current problem, the nonlinearity of the macroscale problem can result from the presence of nonlinear magnetic materials in the microstructure and/or the presence of significant eddy currents that can introduce a phase shift between the magnetic flux density $\bb$ and the magnetic field $\bh$. 

The paper is organized as follows. Section~\ref{sec:finescale_formulations} focuses on the classical $\bh-\phi$ formulation for the finescale problem which will be used as the reference problem. In Section~\ref{sec:dimensional_asymptotic_analysis} we carry out a normalization process for Maxwell's equations allowing for a dimensional analysis of Maxwell's equations and their asymptotic analysis in the context of asymptotic homogenization with lack of scale separation. In Section~\ref{sec:multiscale_formulations}, the strong and weak homogenized formulations are proposed and discretized in space and time leading to a nonlinear system of equations which is linearized using the Newton--Raphson method. The upscalings of the magnetic flux density and of the homogenized incremental permeability are also detailed. Section~\ref{sec:results} deals with numerical examples. At first, we validate the accuracy of the multiscale $\bh-\phi$ formulation by comparing results of the HMM method to the reference solutions obtained by solving the finescale problem. We then study the performance of the method by carrying out strong and weak scaling for the multiscale problems. In Section~\ref{sec:conclusions_perspectives} we close the paper with some conclusions and perspectives. 
%
%
\section{The finescale formulations}
\label{sec:finescale_formulations}
%
%
In this section, we recall the main results of the $\bh-\phi$ formulation for the reference problem. Much of the content of this section can be found in the literature, e.g. in~\cite{hiptmair-mqs-05, dular-formulation-99}.

We consider the magnetoquasistatic problem governed by the following Maxwell equations and constitutive laws~\cite{jackson-em-98}:
\begin{subequations}
  \begin{gather}
      \partial_t \bb^{\varepsilon} + \Curl[]{\be^{\varepsilon}} = \bzero , \,\,\,\, \Curl[]{\bh^{\varepsilon}}  = \bj^{\varepsilon} , \,\,\,\, \Div[]{\bb^{\varepsilon}} = 0 \quad \text{ in } \Omega, 
      \label{eq:maxwell-equations}
      \\
      \bj^{\varepsilon} = \sigma^{\varepsilon} \be^{\varepsilon} + \bj_s , \,\,\,\, \bb^{\varepsilon}=\bBB\big( \bh^{\varepsilon} \big) = \mu(\bh^{\varepsilon}) \bh^{\varepsilon} \quad \text{ in } \Omega, 
      \label{eq:constitutive-laws}
        \\
        \bn \times \bh^{\varepsilon}\big|_{\Gamma} = \bzero,
        \label{eq:bnd}
        \\
        \bh^{\varepsilon}(\bx, t = 0) = \bh_0^{\varepsilon}(\bx) \quad \text{ in } \Omega.
        \label{eq:ic}
  \end{gather}
  \label{eq:mqs_problem}
\end{subequations}
\begin{figure}
  \centering
    \def\svgwidth{1.0\textwidth}
    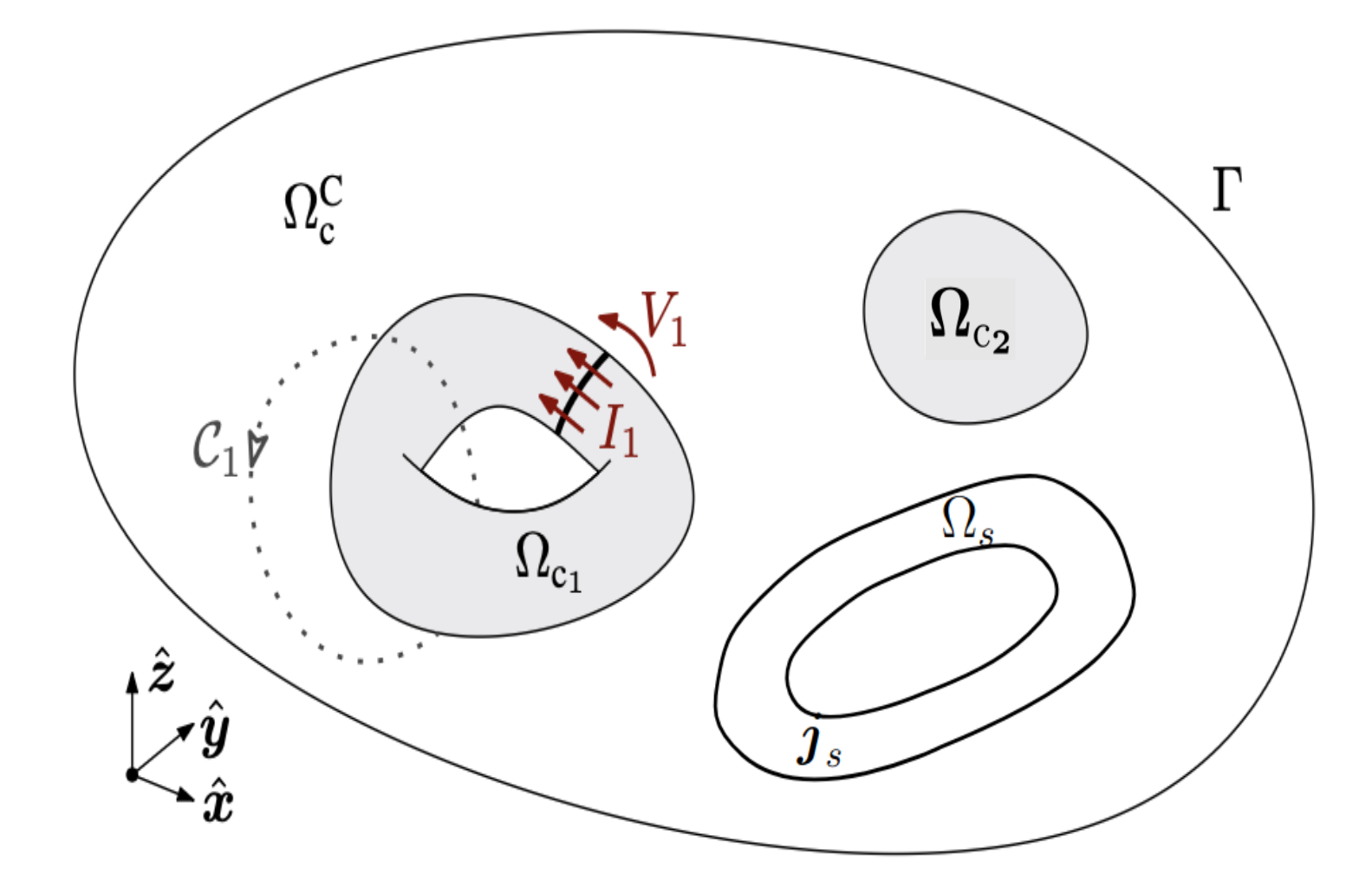
    \caption{Bounded domain $\Omega$ and its subregions~\cite{niyonzima-ah-16, dular-thesis-23}.}
  \label{fig:computational_domains}
\end{figure}
In \eqref{eq:maxwell-equations}--\eqref{eq:ic}, $\bh^{\varepsilon}, \be^{\varepsilon}$, $\bb^{\varepsilon}$ and $\bj^{\varepsilon}$ are the magnetic field [A/m], the electric field [V/m], the magnetic flux density [T] and the current density [A/m$^2$], respectively. The source current density $\bj_s$ is defined in the inductors $\Omega_s$, and $\bn$, $\sigma^{\varepsilon}$ and $\mu^{\varepsilon}$ are the outward normal vector, the electric conductivity [S/m] and the magnetic permeability [H/m], respectively.
The superscript $^{\varepsilon}$ denotes the multiscale nature of the fields. Problem \eqref{eq:maxwell-equations}--\eqref{eq:bnd} is defined on an open bounded domain $\Omega$ of $\mathbb{R}^2$ or $\mathbb{R}^3$, with $\Omega = \Omega_c \cup \Omega_c^C$ and $\Omega_c \cap \Omega_c^C = \emptyset$, where $\Omega_c$ is the conducting domain which contains conducting materials and massive inductors, $\Omega_c^C$ is the non-conducting domain which can be multiply connected (see Figure~\ref{fig:computational_domains}).
A zero tangential trace is prescribed on the boundary $\partial \Omega = \Gamma$ of the domain $\Omega$.

The following assumptions on the data of the problem are made to ensure the existence of a unique solution~\cite{li-h-formulation-13, jiang-homogenization-14, bermudez-math-electromagnetism-14, rodriguez-eddycurrents-10, hiptmair-mqs-05, visintin-tsh-11}:
\begin{enumerate}
  \item[1.] The electric conductivity and the differential magnetic permeability satisfy: 
  \begin{align}
    &0 < \sigma_{\mathrm{min}} \leq \sigma^{\varepsilon}(\bx) \leq \sigma_{\mathrm{max}} \text{ in } \Omega_c \quad \text{ and } \quad \sigma = 0 \text{ in } \Omega_c^C, 
    \\
    &0 < \mu_{\mathrm{min}} \leq \bigg|\ \left( \frac{\partial \bBB}{\partial \bh^{\varepsilon} } \right)_{ij} \bigg| \leq \mu_{\mathrm{max}} \text{ in } \Omega \quad \text{ for } i, j = 1, 2, 3.
  \end{align}
  \item[2.] 
  Strong monotonicity of $\bBB$, i.e., $\exists C_1 >0$ such that: 
  \begin{equation}
    (\bBB(\bv_2) - \bBB(\bv_1))\cdot(\bv_2 - \bv_1) \geq C_1 | \bv_2 - \bv_1 |^2 \,\, \forall \bv_1, \bv_2 \in \mathbb{R}^3.
  \end{equation}
  \item[3.] 
  Lipschitz continuity of $\bBB$, i.e., $\exists C_2 >0$ such that:
  \begin{equation}
    |\bBB(\bv_2) - \bBB(\bv_1)| \leq C_2 | \bv_2 - \bv_1 | \,\, \forall \bv_1, \bv_2 \in \mathbb{R}^3.
  \end{equation}
  \item[4] For almost every instant $t \in [0, T]$, $\bj_s(\cdot, t) \in \LLtwo[]{\Omega}$ and $\Div[]{\bj_s}(\cdot, t) = 0$. These two conditions can be fulfilled for instance if:
  \begin{equation}
    \bj_s(\cdot, t) \in X(\Omega):= \bigg\{ \bu \in \Hdiv[]{\Omega} \,\, \Big| \,\, \Div[]{\bu} = 0 \bigg\}
  \end{equation}
    where the definition of the space $\Hdiv[]{\Omega}$ and other spaces used for electromagnetic fields can be found, e.g., in \cite{monk-fem-03}. 
    As a function of time, we consider sources $\bj_s$ with enough regularity (e.g., $\bj_s(x,\cdot)$ of class $C^1$ which leads to a solution with a $C^1$ regularity in time. This assumption is good e.g., for sine sources, but not good e.g., for PWM sources.
\end{enumerate}
For details regarding the proof of existence, we refer the reader to Appendix A of \cite{li-h-formulation-13} and Section 4 of \cite{visintin-tsh-11}.

The solenoidal condition on the source in Assumption 4 can be numerically achieved by setting $\bj_s(\bx, t) = \Curl[]{\bh_s(\bx, t)} = I(t) \, \Curl[]{\bh_{s0}(\bx)}$ where $I(t)$ is the current flowing through the inductor and $\bh_{s0}$ is the static magnetic field obtained by imposing the current density $\bj_{s0}$ with the total net current $\displaystyle \int_{S} \bj_{s0} \text{d}S$ = 1\,A. This field can be computed by solving a projection problem with the Coulomb or a tree-cotree gauge \cite{creuse-gauge-19}. 

In this paper, we adopt the following notation:
\begin{equation}
  (u, v)_{\Omega}:= \int_{\Omega} u(\bx) \cdot v(\bx)  \text{d}\bx \,\, \text{ and } \,\, <u, v>_{\Gamma}:= \int_{\Gamma} u(\bx) \cdot v(\bx) \text{d}\bx. 
\end{equation}
Combining Maxwell's equations (\ref{eq:mqs_problem}a)--(\ref{eq:mqs_problem}c), we get the following $\bh$-conforming strong form:
\begin{equation}
  \begin{aligned}
    \partial_t \bBB(\bh^{\varepsilon}) + \Curl[]{\left(\rho^{\varepsilon} \Curl[]{\bh^{\varepsilon}} \right)} &= \bzero &&\quad \text{ in } \Omega_c, \\
    \partial_t \Div[]{\bBB(\bh^{\varepsilon})} &= 0 &&\quad \text{ in } \Omega_c^C, \\
    \mathbf{n} \times [\![ \bh^{\varepsilon} ]\!] = \mathbf{0}, \mathbf{n} \cdot [\![ \bb^{\varepsilon} ]\!] &= 0 &&\quad \text{ on } \Sigma_c, \\
    \bn \times \bh^{\varepsilon}| &= \bzero  &&\quad \text{ on } \Gamma, \\
    \bh^{\varepsilon}(\bx, t = 0) &= \bh_0^{\varepsilon}(\bx) &&\quad \text{ in } \Omega,
  \end{aligned}
  \label{eq:strong-eqn}
\end{equation}
The problem can either be fed by known current or voltage \cite{dular-formulation-99, pellikka-homology-13}. For simplicity, we restrict the analysis in this paper to current-driven problems. In \eqref{eq:strong-eqn}, $[\![ \bullet ]\!]$ denotes the jump of a quantity at an interface, and $\Sigma_c$ is the interface between $\Omega_c$ and $\Omega_c^C$. 
The resulting weak form can be written as \cite{hiptmair-mqs-05, dular-formulation-99, dular-thesis-23}: find $\bh^{\varepsilon} \in C^1([0, T]; \bV_I(\Omega) ) $ 
such that 
\begin{equation}
  \Big(\partial_t \bBB(\bh^{\varepsilon}), \bv \Big)_{\Omega} + \Big(\rho \, \Curl[]{\bh^{\varepsilon}}, \Curl[]{\bv} \Big)_{\Omega_{c}} = 0.
  \label{eqref:section_finescale_weak_form_magn}
\end{equation}
holds for all test functions $\bv \in \bV_0(\Omega)$. The spaces $\bV_I(\Omega)$ and $\bV_0(\Omega)$ are defined by 
\begin{equation*}
  \begin{aligned}
    \bV_I(\Omega) &:= \Big\{ \bh \in \Hcurl[]{\Omega} \, \big| \, \Curl[]{\bh} = \bzero \text{ in } \Omega_c^C, \bn \times \bh|_{\Gamma} = \bzero, \mathcal{I}_i(\bh) = I_i \text{ for } i \in C_I \Big\}, \\
    \bV_0(\Omega) &:= \Big\{ \bh \in \Hcurl[]{\Omega} \, \big| \, \Curl[]{\bh} = \bzero \text{ in } \Omega_c^C, \bn \times \bh|_{\Gamma} = \bzero , \mathcal{I}_i(\bh) = 0 \text{ for } i \in C_I \Big\},
  \end{aligned}
\end{equation*}

\noindent where the functional $\mathcal{I}_i(\bh)$ denotes the net current $I_i$ flowing in the subsets $\Omega_{ci} \subseteq \Omega_c$, $i = 1, \dots, N_{\mathrm{mass}}$ and $\Omega_{si} \subseteq \Omega_s$, $i = 1, \dots, N_{\mathrm{str}}$ where a current is imposed. This functional is given as the circulation of $\bh$ along a closed loop $\mathcal{C}_i$ around that conductor:
\begin{equation*}
  \mathcal{I}_i(\bh) = \oint_{\mathcal{C}_i} \bh \cdot d\ell = I_i.
\end{equation*}
The set $C_I$ consists of the curves $\mathcal{C}_i$ for $i = 1, 2, \dots, N_{\mathrm{mass}} + N_{\mathrm{str}}$. 
The space of solutions $C^1([0, T]; \bV_I(\Omega))$ is defined by (see section 5.9. of~\cite{evans-pde-10} or chap.~23 of~\cite{zeidler-nfa-13-1}): 
\begin{multline*}
  C^1([0, T]; \bV_I(\Omega)):= \Bigg\{\bu: [0, T] \to \bV_I(\Omega): \\ 
  \|\bu\|_{C^1([0, T]; \bV_I(\Omega) ) } = \underset{t \in ]0, T[}{\text{sup}} \Big\|\bu(t) \Big\|_{\bV_I(\Omega) } + \underset{t \in ]0, T[}{\text{sup}} \Big\|\frac{\text{d} \bu(t)}{\text{d} t} \Big\|_{\bV_I(\Omega) } \Bigg\}. 
\end{multline*}
For the characterization of $\bV_I(\Omega)$ see also section~1.2 of \cite{hiptmair-mqs-05, pellikka-homology-13, dular-thesis-23} and references therein.

At the discrete level, the solution $\bh^{\varepsilon}$ can be approximated by~\cite{dular-formulation-99, pellikka-homology-13}:
\begin{multline}
  \bh^{\varepsilon}(\bx, t) \approx \sum_{e \in \Omega_c} h_e(t) \bSS_e(\bx) + \sum_{n \in \Omega_c^C} \phi_n(t) \Grad[]{S_n}(\bx) \\ + \sum_{i = 1}^{N_{\mathrm{mass}}} I_i(t) \bc_i(\bx) + \sum_{i = 1}^{N_{\mathrm{str}}} I_i(t) \bh_{s0, i}(\bx),
  \label{eq:decomposition_bh}
\end{multline}
where the first term of \eqref{eq:decomposition_bh} only defined in conductors $\Omega_c$ is approximated using first order Nédélec elements defined on an element $T$ by:
$$
  \mathcal{ND}_1(T) := \left\{ \bFF_h : T \to \mathbb{R}^3, \bx \mapsto \ba + \bb \times \bx, \ba, \bb \in \mathbb{R}^3 \right\}.
$$
The second term only defined in non-conducting domains $\Omega_c^C$ is approximated as a gradient of a nodal function using classical $\mathbb{P}_1(\mathcal{T}_h)$ Lagrange elements \cite{ern-fem-21}. 
The third term accounts for the topology of $\Omega_c^C$ via global basis functions associated with \emph{cuts} that render the domain simply connected \cite{dular-formulation-99}. This is defined for each of the $N_{\mathrm{mass}}$ massive inductors. Finally, the last term represents the contribution to the magnetic field of the $N_{\mathrm{str}}$ stranded inductors. 
%
%
\section{Dimensional and asymptotic analysis}
\label{sec:dimensional_asymptotic_analysis}
%
%
The challenges of deriving the homogenized model for eddy current problems are twofold: (1) the violation of the so-called \emph{scale separation} assumption; (2) the derivation of the governing PDEs of the macroscale and cell problems that can be of different nature in the case of the locally-confined eddy current problem. In this section, we conduct an asymptotic analysis of the normalized Maxwell equations for eddy current problems to show, that in the case of the eddy current multiscale problem, the \emph{scale separation} assumption implicitly assumed for all classical asymptotic homogenization methods can be violated if the skin depth is comparable to the size of the cell. We also take a preliminary look at the second challenge which will be further treated in Section~\ref{sec:multiscale_formulations} by proposing a well adapted homogenized model for locally-confined eddy current problems.
%
%
\subsection{Dimensional and asymptotic analysis}
\label{sec:dimensional_asymptotic_analysis_1}
%
%
The \emph{scale separation} assumption implies that all the characteristic lenghts of the physical problem must be very large compared to the size of the periodic cell $l_{c, i}$ along the $i^{\mathrm{th}}$ direction with $i = 1, 2, 3$ for periodic homogenization or the size of the representative volume element (RVE) for stochastic homogenization. In the context of Maxwell's equations, many such characteristic lengths can be defined. They include the macroscale characteristic length of the device $L_{c, i}$ along the $i^{\mathrm{th}}$ direction with $i = 1, 2, 3$ for 3D problems, the wavelength of the exciting source $\lambda = 1/(f \sqrt{\mu \epsilon})$ and the skin depth $\delta$. The justification of the eddy current problem from Maxwell's equations has already been conducted in~\cite{rodriguez-eddycurrents-10, buffa-eddycurrents-00, rapetti-mqs-14}. Here we illustrate it in the context of the eddy current multiscale problem by applying a normalization process similar to the one used in Section 4.2.2. of~\cite{niyonzima_multiphysics-19} and along the lines of \cite{amirat-homogenizationMaxwell-11, amirat-homogenizationMaxwell-17, bouvet-multiscale-23} to the following full linear Maxwell equation: 
\begin{equation}
  \epsilon \frac{\partial^2 \be}{\partial t^2} + \sigma \frac{\partial \be}{\partial t} + \Curl[]{ \left( \frac{1}{\mu} \Curl[]{\be} \right) } = -\frac{\partial \bj_s}{\partial t} \quad \text{ in } \Omega.
  \label{eq:maxwell_normalization_e}
\end{equation}
In the rest of the section, quantities with an overbar are dimensionless while quantities with the index $_c$ are characteristic quantities.
Defining $\be = E_c \bar{\be}$, where $E_c$ is the characteristic value of the electric field [V/m] and $\bar{\be}$ is the dimensionless electric field, and introducing a new dimensionless coordinate system ($\tau$, $\boldsymbol{\eta}$) in place of ($t$, $\bx$) where $t = T_{c} \, \tau$ and $x_i = L_{c_{i}} \, \eta_i$ for $i = 1, 2$ and $3$, we get:
\begin{equation}
  \begin{aligned}
      &dt = T_{c} \, d \tau \, \, , \, \,  
      \frac{\partial (\cdot)}{\partial t} = \frac{1}{T_{c}} \frac{\partial (\cdot)}{\partial \tau} \, \, , \\ 
      &dx_i = L_{c_{i}} \, d \eta_i \, \, , \, \, 
      \displaystyle \frac{\partial (\cdot)}{\partial x_i} = \frac{1}{L_{c, i}} \frac{\partial (\cdot)}{\partial \eta_i} \, \, , \, \, 
      \displaystyle \frac{\partial^2 (\cdot)}{\partial x_i^2} = \frac{1}{L_{c, i}^2} \frac{\partial^2 (\cdot)}{\partial \eta_i^2}, \\
      &\displaystyle \Curl[]{( \frac{1}{\mu} \Curl[]{\be} ) } = \frac{E_c}{\mu_c L_c^2} \Curl[_{\bbeta}]{ \left( \frac{1}{\bar{\mu}} \Curl[_{\bbeta} ] {\bar{\be} } \right) },
    \label{eq:maxwell_normalization_2}
  \end{aligned}
\end{equation}
with $L_c = \underset{i = 1, 2, 3}{\text{min}} (L_{c_{i}})$. This leads to the following normalized equation: 
\begin{equation}
  \label{eq:maxwell_adim_1}
    \frac{ (\epsilon_c \bar{\epsilon}) }{T_c^2} \frac{\partial^2 \bar{\be}}{\partial \tau^2} 
    + \frac{(\sigma_c \bar{\sigma})}{T_c} \frac{\partial \bar{\be}}{\partial \tau}
    + \frac{1}{\mu_c L_c^2} \Curl[_{\bbeta} ]{ \left( \frac{1}{\bar{\mu} } \Curl[_{\bbeta} ] {\bar{\be} } \right) } 
     = - \frac{1}{T_c} \frac{J_c}{E_c} \frac{\partial \bar{\bj}_s}{\partial \tau}.
\end{equation}
The parameters $\tau$ and $\eta_i$, $i$ = 1, 2 and 3, are dimensionless temporal and spatial coordinates. The notation $\Curl[_{\bbeta}]{(\cdot)}$ is used for the curl with respect to the dimensionless coordinates $\bbeta$. 
All quantities with the subscript $_c$ are the characteristic quantities for the corresponding physical quantities, i.e., $E_c, J_c, T_{c}, \mu_c, \epsilon_c$ and $\sigma_c$ are characteristic electric field, electric current density, time, permeability, permittivity and conductivity, respectively. The characteristic time $T_c$ can be considered as the period of the exciting source which is related to the charasteristic fundamental frequency by $f_c = 1/T_c = \omega_c/(2 \, \pi)$ where $\omega_c$ is angular velocity. Defining the characteristic velocity of the electromagnetic waves as $v_c = 1/\sqrt{(\mu_c \epsilon_c)}$, the characteristic wavelength is given by $\lambda_c = v_c/f_c = 2 \, \pi \, v_c/\omega_c$ and the characteristic skin depth by: $\delta_c = \sqrt{2/(\sigma_c \omega_c \mu_c)}$.
This leads to the following normalized equation:
\begin{equation}
  \begin{aligned}
      & \frac{ (\epsilon_c \bar{\epsilon}) }{T_c^2} \frac{\partial^2 \bar{\be}}{\partial \tau^2} 
    + \frac{(\sigma_c \bar{\sigma})}{T_c} \frac{\partial \bar{\be}}{\partial \tau}
    + \frac{1}{\mu_c L_c^2} \Curl[_{\boldsymbol{\eta}}]{ \left( \frac{1}{\bar{\mu} } \Curl[_{\boldsymbol{\eta} } ] {\bar{\be} } \right) } \underset{\displaystyle [\times \mu_c]}{=}
      \\
      & \frac{\epsilon_c \mu_c}{T_c^2} \bar{\epsilon} \frac{\partial^2 \bar{\be}}{\partial \tau^2} 
      + \frac{\mu_c \sigma_c}{T_c} \bar{\sigma} \frac{\partial \bar{\be}}{\partial \tau}
      + \frac{1}{L_c^2} \Curl[_{\boldsymbol{\eta}}]{ \left( \frac{1}{\bar{\mu} } \Curl[_{\boldsymbol{\eta}}]{\bar{\be} } \right) } =
      \\    
      & \frac{\omega_c^2}{4 \, \pi^2 \, v_c^2} \bar{\epsilon} \frac{\partial^2 \bar{\be}}{\partial \tau^2} 
      + \frac{\mu_c \sigma_c}{T_c} \bar{\sigma} \frac{\partial \bar{\be}}{\partial \tau}
      + \frac{1}{L_c^2} \Curl[_{\boldsymbol{\eta}}]{ \left( \frac{1}{\bar{\mu} } \Curl[_{\boldsymbol{\eta}}]{\bar{\be} } \right) } =
      \\    
      & \frac{1}{\lambda_c^2} \bar{\epsilon} \frac{\partial^2 \bar{\be}}{\partial \tau^2} 
      + \frac{1}{\pi \, \delta_c^2} \bar{\sigma} \frac{\partial \bar{\be}}{\partial \tau}
      + \frac{1}{L_c^2} \Curl[_{\boldsymbol{\eta}}]{ \left( \frac{1}{\bar{\mu} } \Curl[_{\boldsymbol{\eta}}]{\bar{\be} } \right) } 
      = - \frac{\mu_c}{T_c} \frac{J_c}{E_c} \frac{\bar{\partial \bj_s}}{\partial \tau }
      \label{eq:maxwell_normalization_2_1}
  \end{aligned}
\end{equation}
where the three terms of the left-hand side represent electric displacement phenomena, eddy currents and magnetic phenomena, respectively. 

It is well known that the ratios of characteristic lengths $\lambda_c/L_c$ and $\delta_c/L_c$ determine the weigh of displacement currents and eddy currents terms. But, they may sometimes not inform on how the local electromagnetic couplings impact the homogenized magneitic behaviour, as will be shown later. 

From now on, we focus on eddy current problems and neglect displacement currents. This assumption is valid if:
\begin{equation}
  \frac{\lambda_c}{L_c} \gg 1 \quad \text{ and } \quad \frac{\epsilon_c \omega_c}{\sigma_c} \ll 1,
\end{equation} 
where the second inequality takes into account the relaxation of electric charges and is obtained from the normalization of $\Div[]{(\bj + \partial_t \bd)} = 0$ and by neglecting the first term of \eqref{eq:maxwell_normalization_2_1}. In this context, this leads to the following normalized equation for multiscale the eddy current problem:
\begin{equation}
    \frac{1}{\pi \, \delta_c^2} \bar{\sigma} \frac{\partial \bar{\be}^{\varepsilon} }{\partial \tau}
    + \frac{1}{L_c^2} \Curl[_{\boldsymbol{\eta}}]{ \left( \frac{1}{\bar{\mu} } \Curl[_{\boldsymbol{\eta}}]{\bar{\be}^{\varepsilon} } \right) }
    = - \frac{\mu_c}{T_c} \frac{J_c}{E_c} \frac{\bar{\partial \bj_s}}{\partial \tau}.
    \label{eq:maxwell_normalization_2_bis}
\end{equation}

The topology of the non-conducting subdomain $\Omega_c^{C \, \varepsilon}$ of the multiscale domain $\Omega^{\varepsilon} \subset \Omega$ plays a key role in the homogenization process of equation \eqref{eq:maxwell_normalization_2_bis} and in determining the type of the PDEs to be solved at the macroscale and the mesoscale levels (see Table 1 and Fig. 1 of~\cite{niyonzima-chm-18}). Two categories of multiscale eddy current problems can be distinguished:
\begin{enumerate}
  \item Eddy current problems with a fully-conducting multiscale domain $\Omega_c^{\varepsilon} \equiv \Omega^{\varepsilon}$, or a multiply-connected conducting domain $\Omega_c^{\varepsilon} \subsetneq \Omega^{\varepsilon}$, which leads to a homogenized problem with macroscale eddy currents $\bj_M$. These problems can be homogenized by upscaling either the material properties such as $\bsigma_M$ and $\bmu_M^{-1}$ or by homogenizing electromagnetic fields such as the magnetic field $\bh_M = \bHH_M(\bb_M, \be_M)$ and the eddy currents $\bj_M = \bJJ_M(\bb_M, \be_M)$ for the $\bh$-conforming formulations, or the magnetic flux density $\bb_M = \bBB_M(\bh_M, \bj_M)$ and the electric field $\be_M = \bEE_M(\bh_M, \bj_M)$ for the $\bb$-conforming formulations.
  \item Locally-confined eddy current problems with a disconnected conducting domain $\Omega_c^{\varepsilon} \subsetneq \Omega^{\varepsilon}$ which leads to a magnetostatic macroscale problem. For this class of problems, we advise upscaling techniques of the electromagnetic fields such as the magnetic field $\bh_M = \bHH_M(\bb_M)$ or the magnetic flux density $\bb_M = \bBB_M(\bh_M)$ instead of the material properties.
\end{enumerate}
Problems of category 2 are the most difficult to deal with as their homogenization leads to a homogenized problem with a magnetostatic macroscale problem and magnetoquasistatic mesoscale problems. 

Two extreme cases can also be distinguished for problems of category 1 with a similar behaviour to that of problems of category 2:
\begin{enumerate}[label=\alph*.]
  \item Problems with global currents and a macroscale skin depth $\delta_M$ comparable to/or even smaller than the cell size $l_{c, i}$, $i$ = 1, 2, 3. The homogenization of these problems makes little sense.
  \item Problems with global currents and a high contrast of conductivities in the matrix and the inclusions, i.e., $\sigma_{\text{matrix}}/\sigma_{\text{inclusion}} \ll 1$. These problems show some features of problems of category 2 when the ratio $\sigma_{\text{matrix}}/\sigma_{\text{inclusion}}$ tends to zero. However, similar problems, e.g., in heat conduction have been treated using some advanced asymptotic expansion methods (see section 4.2 of \cite{auriault-homogenization-10}). 
\end{enumerate}
While it is possible to derive a homogenized problem for problems of category 1 using the classical asymptotic expansion method, it is not the case for problems of category 2 due to the vanishing homogenized conductivity. In the following paragraph we briefly introduce the asymptotic expansion method well adapted for problems of category 1 and in Section~\ref{sec:multiscale_formulations} we propose a multiscale model able to handle problems of category 2.

In a classical asymptotic expansion method with two separated scales which naturally leads to the homogenization of the material property (see e.g.~\cite{fish-homogenization-13, bensoussan-ahm-78, sanchez-ahm-80, zhikov-conv-91, allaire-tsh-92}), a distinction is made between the macroscale coordinates $\bx = (x_1, x_2, x_3)$ and the mesoscale coordinates $\by = (y_1, y_2, y_3) = \bx/\varepsilon$. The spatial derivatives and derived differential operators such as the $\GradSymb$ and the $\CurlSymb$ can then be expanded using:  
\begin{equation}
  \frac{\partial (\cdot)^{\varepsilon}}{\partial x_i} = 
  \frac{\partial (\cdot)}{\partial x_i} + \frac{\partial y_i}{\partial x_i} \frac{\partial (\cdot)}{\partial y_i}
  = \frac{\partial (\cdot)}{\partial x_i} + \frac{1}{\varepsilon_i} \frac{\partial (\cdot)}{\partial y_i},
    \label{eq:asymptotic_expansion_1}
\end{equation}  
where $\varepsilon_i = l_{c, i}/L_{c, i}$ is the ration between the mesoscale and the macroscale characteristic lengths along the direction $i$. 
Assuming a multiscale problem of category 1, the unknown electric field in~\eqref{eq:maxwell_normalization_2_bis} can be expanded in terms of powers of $\varepsilon$ as:
\begin{multline}
  \bar{\be}^{\varepsilon}(\bx, t(\tau)) = \sum_{i = 0}^{\infty} \varepsilon^i \bar{\be}_i (\bbeta_{\bx}, \bbeta_{\by}, \tau) = 
  \\
  \bar{\be}_0 (\bbeta_{\bx}, \bbeta_{\by}) + \varepsilon \bar{\be}_1 (\bbeta_{\bx}, \bbeta_{\by}, \tau) + \varepsilon^2 \bar{\be}_2 (\bbeta_{\bx}, \bbeta_{\by}, \tau) + \dots
  \label{eq:e_expansion}
\end{multline}
where $\bbeta_{\bx} := \bx/L_c$ and $\bbeta_{\by} := \bx/l_c = \bx/(\varepsilon L_c) = \by/L_c$ are the dimensionless macroscale and mesoscale coordinates. 
Using~\eqref{eq:asymptotic_expansion_1} in~\eqref{eq:maxwell_normalization_2_bis} and these definitions, we get:
\begin{equation}
  \begin{aligned}
      \frac{1}{\pi \, \delta_c^2} \bar{\sigma} &\frac{\partial \bar{\be}^{\varepsilon}}{\partial \tau}(\bx, t(\tau)) 
      + \frac{1}{L_c^2} \Curl[_{\boldsymbol{\eta}}]{ \left( \frac{1}{\bar{\mu}} \Curl[_{\boldsymbol{\eta}}]{\bar{\be}^{\varepsilon} ) } \right) } 
       = \frac{1}{\pi \, \delta_c^2} \bar{\sigma}(\by) \frac{\partial }{\partial \tau} \left( \bar{\be}^{\varepsilon} \right) 
       \\
       &+ \frac{1}{L_{c}^2} \Curl[_{\bbeta_{\bx}}]{ \left( \frac{1}{\bar{\mu}(\by) } \Curl[_{\bbeta_{\bx}}]{ \bar{\be}^{\varepsilon} }\right) }
      + \frac{1}{\varepsilon L_{c}^2} \Curl[_{\bbeta_{\bx}}]{ \left( \frac{1}{\bar{\mu}(\by) } \Curl[_{\bbeta_{\by}}]{ \bar{\be}^{\varepsilon}  } \right) }
      \\
       &+ \frac{1}{\varepsilon L_{c}^2} \Curl[_{\bbeta_{\by}}]{ \left( \frac{1}{\bar{\mu}(\by) } \Curl[_{\bbeta_{\bx}}]{ \bar{\be}^{\varepsilon}  }\right)}
      + \frac{1}{\varepsilon^2 L_{c}^2} \Curl[_{\bbeta_{\by}}]{ \left( \frac{1}{\bar{\mu}(\by) } \Curl[_{\bbeta_{\by}}]{ \bar{\be}^{\varepsilon} } \right) }
      \\= &- \frac{\mu_c}{T_c} \frac{J_c}{E_c} \frac{\bar{\partial \bj_s}}{\partial \tau },
      \label{eq:maxwell_normalization_3}
  \end{aligned}
\end{equation}
where $\bar{\be}^{\varepsilon}=\bar{\be}^{\varepsilon}(\bx, t(\tau))$, which can further be expanded in terms of powers of $\varepsilon$ using~\eqref{eq:e_expansion}. In~\eqref{eq:maxwell_normalization_3}, $\Curl[_{\bbeta_{\bx}}]{}$ and $\Curl[_{\bbeta_{\by}}]{}$ are curl operators applied to the macroscale and mesoscale dimensionless coordinates, respectively. 
In most textbooks dealing with the asymptotic expansion of parabolic equations, the term corresponding to the time derivative is often considered to be of order $\mathcal{O}(\varepsilon^0)$. The equation corresponding to order $\mathcal{O}(\varepsilon^{-2})$ is therefore used to justify the fact that the first term of the expansion (the macroscale field) is independent from the mesoscale coordinate system $\bar{\be}_0(\bbeta_{\bx}, \bbeta_{\by}, \tau) = \bar{\be}_0(\bbeta_{\bx}, \tau)$, while the equations of orders $\mathcal{O}(\varepsilon^{-1})$ and $\mathcal{O}(\varepsilon^{0})$ are respectively used to define the cell problem and the macroscopic problem. These developments assume that the following scale separation hypotheses are satisfied
\begin{equation}
  \!\!\!\! \varepsilon \ll 1, \,\,\,\, l_{c, i} = \varepsilon L_{c, i} \ll \delta_c \,\,\,\, \text{ and } \,\,\,\, l_{c, i}^{i} = \varepsilon L_{c, i} \ll \lambda_c \,\,\,\, \text{ for } i = 1, 2, 3,
  \label{eq:maxwell_characteristic_lengths}
\end{equation}
which implies that the characteristic macroscale length $L_{c}$, the skin depth $\delta_c$ and the wave length $\lambda_c$ are very large compared the mesoscale characteristic length $l_{c}$. Unfortunately the second inequality of \eqref{eq:maxwell_characteristic_lengths} may not be verified for eddy current problems as $l_{c} = \varepsilon L_{c} \approx \delta$. Therefore, a rigorous asymptotic expansion for eddy current problems without scale separation should be carried along the lines developed in \cite{amirat-homogenizationMaxwell-11, amirat-homogenizationMaxwell-17} for time harmonic problems.
\begin{figure}[htbp]
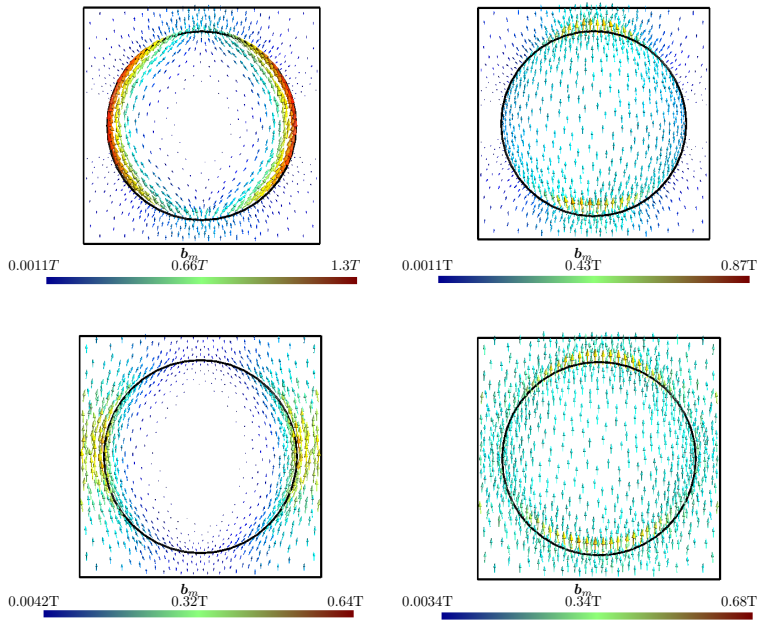

    \scalebox{0.55}{
        \hspace{30mm} 
        \input{./data/b_f1e4_mur1e4_Other_TS15_New.tex} \hspace{-60mm}
        \input{./data/b_f1e-2_mur1e4_Other_TS15_New.tex} \hspace{-60mm}
        }
        
      \scalebox{0.55}{
        \hspace{30mm} 
        \input{./data/b_f1e8_mur1e0_Other_TS15_New.tex} \hspace{-60mm}
        \input{./data/b_f1e2_mur1e0_Other_TS15_New.tex} \hspace{-60mm}
        }
    \caption{\footnotesize
    Maps of the magnetic flux density of magnetoquasistatic and magnetostatic problems with a conducting magnetic inclusion of radius $r = 40 \mu$m. \textbf{Top left}: the magnetic flux density $\bb_m^{\mathrm{MQS}}$ at $t = 4.6875 \times 10^{-6}$ s obtained solving a magnetoquasistatic problem with $f = 10^4$ Hz, $\mu_r = 10^4$ and $\delta = 5\mu$m. \textbf{Top right}: the magnetic flux density $\bb_m^{\mathrm{MS}}$ obtained solving a magnetostatic problem with $\mu_r = 10^4$. \textbf{Bottom left}: the magnetic flux density $\bb_m^{\mathrm{MQS}}$ at $t = 4.6875 \times 10^{-10}$ s obtained solving a magnetoquasistatic problem with $f = 10^8$ Hz, $\mu_r = 1$ and $\delta = 5\mu$m. \textbf{Bottom right}: the magnetic flux density $\bb_m^{\mathrm{MS}}$ obtained solving a magnetostatic problem with $\mu_r = 1$.
    The magnetoquasistatic problem in Figure (top-left) and the magnetostatic problem in Figure (top-right) yield the same homogenized magnetic flux density (see Figure \ref{fig:section3_img2} (top-right)) whereas the magnetoquasistatic problem in Figure (bottom-left) and the magnetostatic problem in Figure (top-right) yield different homogenized magnetic flux density (see Figure \ref{fig:section3_img2} (bottom-right). This seems to corroborate the influence of the ratio $\mathrm{P}_{m}^{\mathrm{eddy}}/\mathrm{P}_{m}^{\mathrm{mag}}$ on the equality of homogenized quantities by the static and the dynamic problems.)
    }\label{fig:section3_img1}
\end{figure}
\normalsize
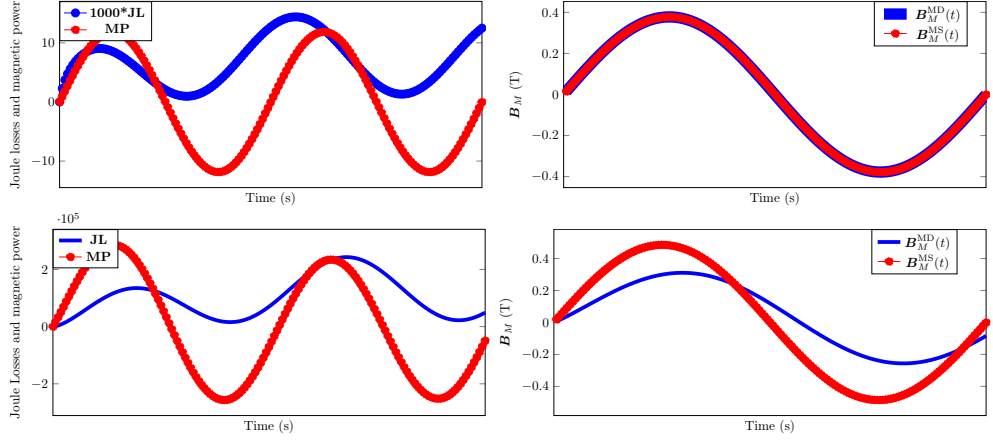
\begin{figure}[htbp]
  \begin{tikzpicture}[scale=0.5]
    \begin{axis}[xlabel={Time (s)}, ylabel={Joule losses and magnetic power}, xmin=0.0, xmax=0.0001, xtick={-1e-3, 1e-3}, width=0.98\columnwidth, height=0.5\columnwidth, legend style={at={(0,1)},anchor=north west}]
      \addplot [color=blue, mark=*, mark size=3.0, mark options=solid] table [x expr=\thisrowno{0}*1, y expr=\thisrowno{1}*1000, col sep=space] {./data/JouleLosses_f1e4_mur1e4_New.csv};
      \addlegendentry{\textbf{1000*JL}}
      \addplot [color=red, mark=*, mark size=3.0, mark options=dashdotted] table [x expr=\thisrowno{0}*1, y expr=\thisrowno{1}, col sep=space] {./data/MagPower_f1e4_mur1e4_New.csv};
      \addlegendentry{\textbf{MP}}
    \end{axis}
  \end{tikzpicture}
  \begin{tikzpicture}[scale=0.5]
    \begin{axis}[xlabel={Time (s)}, ylabel={$\bBB_M$ (T)}, xmin=0.0, xmax=100, xtick={-1e3, 1e3}, width=0.98\columnwidth, height=0.5\columnwidth, legend style={at={(0.735,1)},anchor=north west}]
      \addplot [color=blue, mark=none, line width=3.0mm, mark size=2.0, mark options=solid] table [x expr=\thisrowno{0}*1e6, y expr=\thisrowno{2}, col sep=space] {./data/b_homog_f1e4_mur1e4_New.csv};
      \addlegendentry{$\bBB_M^{\mathrm{MQS}}(t)$}
      \addplot [color=red, mark=*, mark size=3.0, mark options=dashdotted] table [x expr=\thisrowno{0}*1, y expr=\thisrowno{2}, col sep=space] {./data/b_homog_f1e-2_mur1e4_New.csv};
      \addlegendentry{$\bBB_M^{\mathrm{MS}}(t)$}
    \end{axis}
  \end{tikzpicture}
  \begin{tikzpicture}[scale=0.5]
    \begin{axis}[xlabel={Time (s)}, ylabel={Joule Losses and magnetic power}, xmin=0.0, xmax=1e-8, xtick={-1e-7, 1e-7},  width=1.0\columnwidth, height=0.5\columnwidth, legend style={at={(0,1)},anchor=north west}]
      \addplot [color=blue, mark=none, line width=1.0mm, mark size=2.0, mark options=solid] table [x expr=\thisrowno{0}*1, y expr=\thisrowno{1}*1, col sep=space] {./data/JouleLosses_f1e8_mur1e0_New.csv};
      \addlegendentry{\textbf{JL}}
      \addplot [color=red, mark=*, mark size=3.0, mark options=dashdotted] table [x expr=\thisrowno{0}*1, y expr=\thisrowno{1}, col sep=space] {./data/MagPower_f1e8_mur1e0_New.csv};
      \addlegendentry{\textbf{MP}}
    \end{axis}
  \end{tikzpicture}
  \begin{tikzpicture}[scale=0.5]
    \begin{axis}[xlabel={Time (s)}, ylabel={$\bBB_M$ (T)}, xmin=0.0, xmax=1e-2, xtick={-1e-1, 1e-1}, width=1.0\columnwidth, height=0.5\columnwidth, legend style={at={(0.735,1)},anchor=north west}]
      \addplot [color=blue, mark=none, line width=1.0mm, mark size=2.0, mark options=solid] table [x expr=\thisrowno{0}*1e6, y expr=\thisrowno{2}, col sep=space] {./data/b_homog_f1e8_mur1e0_New.csv};
      \addlegendentry{$\bBB_M^{\mathrm{MQS}}(t)$}
      \addplot [color=red, mark=*, mark size=3.0, mark options=dashdotted] table [x expr=\thisrowno{0}*1, y expr=\thisrowno{2}, col sep=space] {./data/b_homog_f1e2_mur1e0_New.csv};
      \addlegendentry{$\bBB_M^{\mathrm{MS}}(t)$}
    \end{axis}
  \end{tikzpicture}
  \caption{\footnotesize
  Influence of the ratio $\mathrm{P}_{m}^{\mathrm{eddy}}/\mathrm{P}_{m}^{\mathrm{mag}}$ on the homogenized magnetic flux density $\bBB_M$ for a conducting magnetic inclusion of radius $r = 40 \mu$m. \textbf{Top left}: Joule losses amplified by a factor $10^3$ and the magnetic power obtained solving a magnetoquasistatic problem with $f = 10^4$ Hz, $\mu_r = 10^4$ and $\delta = 5\mu$m. 
  \textbf{Top right}: The homogenized magnetic flux density $\bBB_M$ upscaled from solutions of the magnetostatic and magnetoquasistatic problems. 
  \textbf{Bottom left}: Joule losses and the magnetic power obtained solving a magnetoquasistatic problem with $f = 10^8$ Hz, $\mu_r = 1$ and $\delta = 5\mu$m. \textbf{Bottom right}: The homogenized magnetic flux density $\bBB_M$ upscaled from solutions of the magnetostatic and magnetoquasistatic problems. 
  }
  \label{fig:section3_img2}
\end{figure}
\normalsize
%
%
\subsection{Upscaling of the magnetic flux density and the powers in the cell}
\label{sec:dimensional_asymptotic_analysis_2}
%
%
In this section, we explain why the influence of dynamic effects on the macroscopic magnetic law is not necessarily related to the skin depth $\delta$ in the conductive inclusion, but rather on the ratio between the eddy currents and the magnetic power in the cell.

In the context of magnetoquasistatic problems, volume average of $\bb_m$ and $\bj_m$ can be used to define the homogenized quantities $\bBB_M$ and $\bJJ_M$, whereas the volume averages of $\bh_m$ and $\be_m$ can provide wrong values of homogenized quantities $\bHH_M$ and $\bEE_M$ \cite{meunier-chm-10}. Accurate upscaling techniques were recently proposed to extend the use of homogenization to the $\ba-v$ formulation that require the upscaling of the magnetic field $\bHH_M$. These techniques can easily be extended to formulations that use the homogenized electric field $\bEE_M$ at the macroscale. 

In  \cite{marteau-hmm-23, marteau-phdthesis-24}, the authors proposed to solve the insulated cell problem and used either edge or face averaging, approaches often used in the high frequency electromagnetic community. They also proposed the definitions of problems used to define the magnetization created by eddy currents confined in the periodic cell. The latter can be subtracted from the magnetic field and used to define the macroscale magnetic field $\bHH_M$.

Another approach recently proposed consists in defining the magnetization $\bM$ created by eddy currents the cell centered around $\bx_c$ by: 
\begin{equation} 
  \bM = \frac{1}{|\Omega_{m}|}\int_{\Omega_{m}} \left(\frac{1}{2} (\bx - \bx_c) \times \bj_m\right) \text{d}\Omega_m,
\end{equation} 
with $\Omega_{m}$ the periodic cell, and in removing it from the average of the magnetic field $\bh$ in order to obtain the homogenized magnetic field as \cite{wulfinghoff-hmm-24}:
\begin{equation}
  \bHH_M = \frac{1}{|\Omega_{m}|}\int_{\Omega_{m}} \bh_m \text{d}\Omega_m - \bM :=  <\bh_m >_{\Omega_{m}} - \bM,
\end{equation}
where the angle brackets $<\cdot>_{\Omega_{m}}$ are used for the average of a quantity over the cell domain $\Omega_m$.

Within the context of the $\ba-v$ formulation, it has already been observed that the volume average $<\bh_m>_{\Omega_{m}}$ can sometimes yield values close to $\bHH_M$, even in the presence of a strong skin effect within the cell. The same observations seem to suggest that the value of the upscaled magnetic field $\bHH_M$ is less dependent on the ratio between the skin depth and the characteristic length $\delta_c/l_{c}$, but rather on the ratio between the average of the total eddy currents losses $\mathrm{P}_{m}^{\mathrm{eddy}}$ and of the total magnetic power $\mathrm{P}_{m}^{\mathrm{mag}}$ contained in the cell.
These powers are defined by: 
\begin{equation}
  \mathrm{P}_{m}^{\mathrm{eddy}} = \frac{1}{|\Omega_m|} \int_{\Omega_m} (\bj_m \cdot \be_m) \text{d} \Omega 
  \quad \quad \text{ and } \quad \quad 
  \mathrm{P}_{m}^{\mathrm{mag}} = \frac{1}{|\Omega_m|} \int_{\Omega_m} (\bh_m \cdot \partial_t \bb_m) \text{d} \Omega.
  \label{eq:average_meso_quantities}
\end{equation}
In the case of our $\bh$-conforming formulation, analogous conclusions can be drawn replacing the accurate volume average of the magnetic flux density $\bb_m^{\mathrm{MQS}}$ obtained by solving the magnetoquasistatic mesoscale problem, by the volume average of the magnetic flux density $\bb_m^{\mathrm{MS}}$ obtained by solving the magnetostatic mesoscale problem. Our numerical tests suggest that the error resulting from the replacement of the average of $\bb_m^{\mathrm{MQS}}$ by that of $\bb_m^{\mathrm{MS}}$ depends on the ratio $\mathrm{P}_{m}^{\mathrm{eddy}}/\mathrm{P}_{m}^{\mathrm{mag}}$ between the average eddy current losses and the magnetic power rather than on the ratio $\delta_c/L_{c}$ between the skin depth and the characteristic size of the mesoscale cell as it can be seen in Figures~\ref{fig:section3_img1} and~\ref{fig:section3_img2}.

Indeed, the magnetic flux density in Figure~\ref{fig:section3_img1} (top-left) exhibits a small skin depth $\delta = 5 \mu$m compared to the radius of the conducting disk $r = 40\mu$m but the ratio of power $\mathrm{P}_{m}^{\mathrm{eddy}}/\mathrm{P}_{m}^{\mathrm{mag}} = 10^{-3}$ is very small as it can be seen in Figure~\ref{fig:section3_img2} (top-left). Therefore, the homogenized magnetic flux density computed from the magnetoquasistatic problem $\bBB_M^{\mathrm{MQS}}(t)$ and from the magnetostatic problem $\bBB_M^{\mathrm{MS}}(t)$ are in a good agreement as can be seen in the top-right image of Figure~\ref{fig:section3_img2}.
The magnetic flux density in Figure~\ref{fig:section3_img1} (bottom-left) exhibits a small skin depth $\delta = 5 \mu$m but with a ratio of power $\mathrm{P}_{m}^{\mathrm{eddy}}/\mathrm{P}_{m}^{\mathrm{mag}} \approx 1$ which is not negligible as it can be seen in Figure~\ref{fig:section3_img2} (bottom-left). Therefore, the homogenized magnetic flux density computed from the magnetoquasistatic and the magnetostatic problems do not coincide at all as it can also be seen in Figure~\ref{fig:section3_img2} (bottom-right).
Both left images of Figure~\ref{fig:section3_img1} exhibits the same ratios $\delta_c/l_c$, however the ratio of powers seems to be the determining factor for the error between the homogenized inductions from the magnetostatic and the magnetoquasistatic problems. This fact is explained in the following paragraphs. 

In the case of locally-confined eddy current problems and assuming that the total power upscaled from the mesoscale to the macroscale is equal to the average total power in the cell, we get:
\begin{equation}
  \mathrm{P}_{M}^{\mathrm{tot}} = \bh_M \cdot \partial_t \bb_M 
  = 
  \frac{1}{|\Omega_m|} \int_{\Omega_m} \left( \bh_m \cdot \partial_t \bb_m +  \bj_m \cdot \be_m\right) \text{d} \Omega_m = < \mathrm{P}_{m}^{\mathrm{tot}} >_{\Omega_m},
  \label{eq:power_equality}
\end{equation}
where the notation $< \mathrm{P}_{m}^{\mathrm{tot}} >_{\Omega_m}$ is used for the volume average of $\mathrm{P}_{m}^{\mathrm{tot}}$ over $\Omega_m$.
The mesoscale fields can be decomposed in terms of the macroscale and mesoscale corrections as:
\begin{alignat*}{2}
  \bb_m &= \bb_M + \bb_c, \quad &&\be_m = \be_M + \kappa \partial_t \bb_M \times \by + \be_c, \\
  \bj_m &= \xcancel{\bj_M} + \bj_c, \quad &&\bh_m = \bh_M + \underset{\text{confined currents}}{\underbrace{\xcancel{\kappa \bj_M \times \by}}} + \bh_c,
\end{alignat*}
where the macroscale fields $\bb_M, \be_M, \partial_t \bb_M$ and $\bh_M$ are constant on the cell and the correction are periodic. 
If the average eddy current losses are negligible compared to total macroscale power, i.e.: 
\begin{equation}
    <\mathrm{P}_{m}^{\mathrm{eddy}} >_{\Omega_m} = \, < \bj_m \cdot \be_m >_{\Omega_m}  \quad \ll \quad  \bh_M \cdot \partial_t \bb_M = \mathrm{P}_{M}^{\mathrm{tot}},
    \label{eq:powers_assumption}
\end{equation}
then the component of $<\bh_c>_{\Omega_m}$ along the direction of $\partial_t\bb_M$ is negligible compared to that of $\bh_M$ along the same direction, i.e., :
\begin{equation}\label{eq:hc_average_inequality}
  <\bh_c>_{\Omega_m} \cdot\, \partial_t\bb_M\, \ll\, \bh_M \cdot\, \partial_t\bb_M.
\end{equation}
As a consequence, $\bh_M$ can very likely be upscaled by a simple volume average of $\bh_m$ as the magnetization created by the time-dependent source term $\partial_t\bb_M$ in $\bh_c$ is mainly oriented along the direction of $\partial_t\bb_M$. 

To prove \eqref{eq:hc_average_inequality}, we compute:
\begin{equation}
    \begin{aligned}
        \mathrm{P}_{M}^{\mathrm{tot}}
        &\approx \frac{1}{|\Omega_m|} \int_{\Omega_m} \left( (\bh_M + \bh_c) \cdot (\partial_t \bb_M + \partial_t \bb_c) \right) \text{d} \Omega_m \\
        &= \underset{\mathrm{I}}{\underbrace{\bh_M \cdot \partial_t \bb_M}} 
        + \bh_M  \cdot \xcancel{<\partial_t \bb_c>_{\Omega_m} } 
        + \underset{\mathrm{II}}{\underbrace{<\bh_c>_{\Omega_m}\cdot\, \partial_t \bb_M}}
        + \underset{\mathrm{III}}{\underbrace{< \bh_c\cdot \partial_t \bb_c >_{\Omega_m} }},
    \end{aligned}
    \label{eq:equality_powers}
\end{equation}
which implies that $\mathrm{II} = - \mathrm{III}$ using the definition of $\mathrm{P}_{M}^{\mathrm{tot}}$.
By successively applying Maxwell-Faraday, Poynting's theorem and the divergence theorem, the last term of \eqref{eq:equality_powers} can be expanded as:    
\begin{equation}
    \begin{aligned}
        \mathrm{III} 
        &= \frac{1}{|\Omega_m|} \int_{\Omega_m} \left( \bh_c \cdot \partial_t \bb_c \right) \text{d} \Omega_m 
         = -\frac{1}{|\Omega_m|} \int_{\Omega_m} \left( \bh_c \cdot \Curl[_y]{\be_c} \right) \text{d} \Omega_m \\
        &= -\frac{1}{|\Omega_m|} \left( 
         \int_{\Omega_m} \Div[_y]{\left( \bh_c \times \be_c \right)} \text{d} \Omega_m 
        + \int_{\Omega_m} \left( \Curl[_y]{\bh_c} \cdot \be_c \right) \text{d} \Omega_m 
        \right)\\
        &= -\frac{1}{|\Omega_m|} \left( 
         \int_{\partial \Omega_m} \bn \cdot \left( \bh_c \times \be_c \right) \text{d} \gamma_m 
        + \int_{\Omega_m} \left( \bj_c \cdot \be_c \right) \text{d} \Omega_m 
        \right)\\
        &= -\frac{1}{|\Omega_m|} \left( 
        \underset{\mathrm{IV_a}}{\underbrace{\int_{\partial \Omega_m} \be_c \cdot \left( \bn \times \bh_c \right) \text{d} \gamma_m } }
        + \underset{\mathrm{IV_b}}{\underbrace{\int_{\Omega_m} \left( \bj_c \cdot \be_c \right) \text{d} \Omega_m }}
        \right).
    \end{aligned}
    \label{eq:equality_powers_2}
\end{equation}
The first term $\mathrm{IV_a}$ of \eqref{eq:equality_powers_2} is zero due to periodicity of $\bh_c$ and $\be_c$ and the anti-periodicity of $\bn$. Thus,  
\begin{equation*}
    \mathrm{II} = -\mathrm{III} = \frac{1}{|\Omega_m|}\int_{\Omega_m} \left( \bj_c \cdot \be_c \right) \text{d} \Omega_m =\, {\mathrm{P}_{M}^{\mathrm{eddy}}},
\end{equation*}
which leads to $<\bh_c>_{\Omega_m}\cdot\, \partial_t \bb_M \ \ll\ \bh_M \cdot \partial_t \bb_M$ using~\eqref{eq:powers_assumption}. 

In the next section, we propose a multiscale model valid for problems of categories 1 and 2 with linear or nonlinear material laws and eddy currents. 
%
%
\section{The multiscale formulations}
\label{sec:multiscale_formulations}
\begin{figure}[H]
  \begin{center}
    \includegraphics[width=0.55\textwidth]{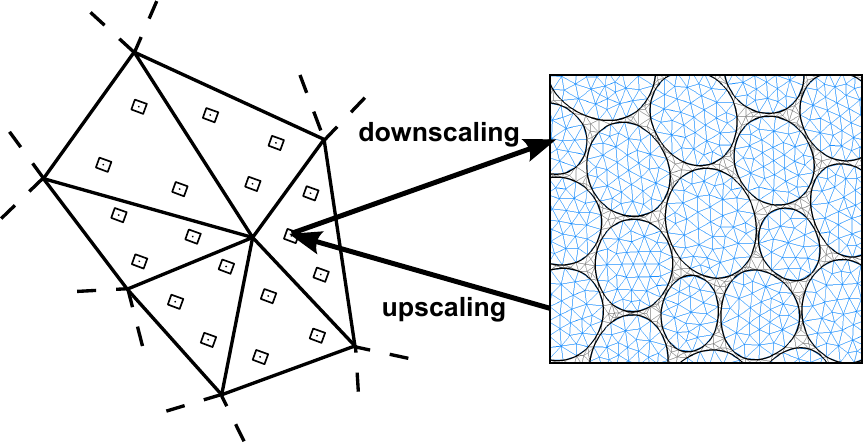}
  \end{center}
  \caption{\footnotesize
  Scale transitions between the macroscale (left) and the mesoscale (right) problems. Downscaling (macro to meso): use of macroscale source terms in the mesoscale problem. Upscaling (meso to macro): computation of homogenized quantities (e.g.\
    material properties or macroscale fields) from the mesoscale solution \cite{niyonzima-hmm-13}. }
  \label{FE2_ppe}
\end{figure}

%
%
In this section we propose an $\bh$-conforming magnetoquasistatic multiscale formulation and its numerical implementation in the context of HMM \cite{e-hmm-03, abdulle-hmm-12, abdulle-hmm-16, niyonzima-hmm-13, niyonzima-chm-14, niyonzima-ah-16}. This method consists in replacing the finescale problem by the homogenized problem which consists of the macroscale problem defined on a coarse mesh together with many mesoscale problems (one per Gau\ss \, point or per barycenter of the macro-element) defined on a fine mesh (see Figure \ref{FE2_ppe}).
\normalsize

The HMM requires a change of subdivisions of the macroscopic domain. The homogenized domain $\Omega_{M,H}$ regroups the heterogeneous materials, that is parts of the former conducting regions $\Omega_c$ (e.g. the conducting inclusions) and of the former non-conducting ones (e.g. insulating material between inclusions). Under the assumption of locally confined eddy currents, the homogenized domain is macroscopically non-conducting, and eddy currents are restricted to the massive inductors. The latter is the remaining non-homogenized conducting regions $\Omega_{M,c} = \Omega_c \backslash \Omega_{M,H}$. Finally, we define $\Omega_{M,a} = \Omega \backslash \left( \Omega_{M,c} \cup \Omega_{M,H} \right)$ the non-homogenized and non-conducting macroscopic regions. Additionally, let us denote the interface between $\Omega_{M,c}$ and $\Omega_{M,a}$ by $\Sigma_{c\, a}$, and that between $\Omega_{M,a}$ and $\Omega_{M,H}$ by $\Sigma_{c\, H}$. For simplicity, it is assumed that $\Sigma_{a\, c}$ and $\Sigma_{c\, H}$ have empty intersection.

Homogenization theories such as the asymptotic expansion method~\cite{bensoussan-ahm-78, fish-homogenization-13}, H-convergence~\cite{murat-gconv-77, tartar-hconv-77}, $\Gamma$-convergence~\cite{DalMaso-gammaconv-93, marcellini-gamma-78} and the two-scale convergence~\cite{nguetseng-tsh-89, allaire-tsh-92, visintin-tsh-06-b, visintin-tsh-08} can be used to derive governing equations at the macro- and meso-scale. For our $\bh$-conforming magnetoquasistatic problems, considering the physics of locally-confined eddy currents, we propose the following homogenized problem for the $\bh$-conforming magnetoquasistatic equations: \\
\begin{subequations}
\noindent The macroscale problem is:
  \begin{equation*}
    \begin{aligned}
        \Curl[_x]{ \left( \rho \, \Curl[_x]{\bh_M} \right) } + \partial_t \bBB_M(\bh_M) &= \bzero \quad \text{ in } \Omega_{M,c},\\
        \Div[_x] \bBB_M(\bh_M)  &= 0 \quad \text{ in } \Omega_{M,a},\\
        \Div[_x] \bBB_M(\bh_M, [\bh_c]^{N^{\mathrm{GP}}})  &= 0 \quad \text{ in } \Omega_{M,H},\\
        \bn \cdot \llbracket\bBB_M(h_M)\rrbracket = 0, \quad \bn \times \llbracket\bh_M\rrbracket &= \bzero \quad \text{ on } \Sigma_{c\, a} \cup \Sigma_{a\, H} \\
        \bn \times \bh_M &= \bzero \quad \text{ on } \Gamma, \\
        \bh_M(\bx, t = 0) &= \bh_0(\bx) \,\, \text{ in } \Omega.
      \end{aligned}
      \label{eq:homogenized-macro-strong}
  \end{equation*}
In $\Omega_{M,c}$ and $\Omega_{M,a}$, the magnetic law $\bBB_M$ is
identical to the original law $\bBB$. In $\Omega_{M,H}$, $\bBB_M$ at
Gau\ss~point $i$ is the volume average of $\bBB(\bh_M + \bh_c^{(i)})$ in the
cell $\Omega_m$, where $\bh_c^{(i)}$ is the solution of the mesoscale problem
at Gau\ss~point $i$ that uses $\bh_M$ as an imput.\\
\noindent The cell problems are:
  \begin{equation*}
    \begin{aligned}
        \partial_t \bBB(\bh_M + \bh_c^{(i)}) + \Curl[_y]{ \left( \rho \Curl[_y]{\bh_c^{(i)}} \right) } &= \bzero \quad \text{ in } \Omega_m,\\
        \bn \times \llbracket\bh_c^{(i)}\rrbracket &= \bzero \quad \text{ on }\Gamma_m \;(\bh_c^{(i)} \text{ is tangentially periodic}),
      \end{aligned}
      \label{eq:homogenized-meso-strong}
  \end{equation*}
\end{subequations}
where $\big[\bh_c\big]^{N^{\mathrm{GP}}}:= \big(\bh_c^{(1)}, \bh_c^{(2)}, \ldots \bh_c^{(N^{\mathrm{GP} } ) } \big)$ is the set of mesoscale solutions computed on periodic cell domains around points of interest of the homogenized domain. In this form, the macroscopic problem is strongly coupled with each mesoscale problem via the circular dependency
$$
  \bh_M \overset{\text{downscale}}{\rightarrow} \bh_M^{\;(i)} \overset{\text{meso. pb.}}{\rightarrow} \bh_c^{(i)} \overset{\text{upscale}}{\rightarrow} \bBB_M \overset{\text{macro. pb.}}{\rightarrow} \bh_M.
$$
\noindent The macroscopic problem will be weakly coupled with the mesoscopic problems once the macroscopic problem is linearized. Linearizing enables solving for the mesoscale problems independently of the macroscopic problem and in parallel, at the cost of an iterative resolution.

The weak form of this problem reads: 
find $(\bh_M, [\bh_c]^{N^{\mathrm{GP}}} ) \in C^1([0, T], \bV_I(\Omega)) \times (C^1([0, T], \bV_{\mathrm{per}}(\Omega_m)))^{N^{\mathrm{GP}}}$ such that for all test functions $(\bv_M, [\bv_c]^{N^{\mathrm{GP}}}) \in \bV_0(\Omega) \times (\bV_{\mathrm{per}}(\Omega_{m}))^{N^{\mathrm{GP}}}$ the following equations hold for $i = 1, 2, \cdots, N^{\mathrm{GP}}$: 
\begin{subequations}
  \begin{align}
    \begin{split}
        &\Big(\partial_t \bBB_M(\bh_M, [\bh_c]^{N^{\mathrm{GP}}} ), \bv_M \Big)_{\Omega} + \Big(\rho \, \Curl[_x]{\bh_M}, \Curl[_x]{\bv_M} \Big)_{\Omega_{c}} = 0,
      \label{eq:homogenized-macro-weak}
    \end{split}
    \\
    \begin{split}
        &\Big(\partial_t \bBB(\bh_M + \bh_{c}^{(i)}), \bv_{c}^{(i)} \Big)_{\Omega_{m}} + \Big(\rho \, \Curl[_y]{\bh_{c}^{(i)}}, \Curl[_y]{\bv_{c}^{(i)}} \Big)_{\Omega_{mc}} = 0.
      \label{eq:homogenized-meso-weak}
    \end{split}
  \end{align}
\end{subequations}

\noindent The space: 
\begin{equation*}
    \bV_{\mathrm{per}}(\Omega_{m}) = \bigg\{ \bu \in \Hcurl[]{\Omega_{m}} \, \Big| \, \Curl{\bu} = \bzero \text{ in } \Omega_{m c}^C ,\,  (\bn \times \bu) \text{ is periodic on } \Gamma_m \bigg\}
\end{equation*}
is the space of $\Hcurl[]{}$-conforming functions on $\Omega_m$ with periodic boundary conditions (PBC). The field $\bh_M$ admits a decomposition similar to \eqref{eq:decomposition_bh}. 
The correction field $\bh_c^{(i)}$ can be approximated by: 
\begin{equation*}
    \bh_c^{(i)}(\by, t) \approx \sum_{e \in \Omega_c^(i)} h_{e, c}^{(i)}(t) \bSS_{e, c}^{(i)}(\by) + \sum_{n \in (\Omega_c^C)^(i)} \phi_{n, c}^{(i)}(t) \Grad[]{S_{n, c}^{(i)}}(\by) 
    + \sum_{i = I}^{N_{\mathrm{mass}}^{\mathrm{m}}} I_{i, c}(t) \bc_{i, c}(\by),
    \label{eq:decomposition_bh_meso}
\end{equation*}
where the last contribution is associated with cuts defined for each of the $N_{\mathrm{mass}}^{m}$ multiply-connected conducting domain in order to render the non-conducting domain $\Omega_{m, i}^{C}$ simply-connected. Regrouping all the macroscale magnetic unknowns from \eqref{eq:decomposition_bh} as:
\begin{equation}
  \bh_M(\bx, t) \approx 
    \underset{\hat{\bh}_M(t)}{\underbrace{
    \begin{pmatrix}
      \hat{h}_{M, 1}(t) & \hat{h}_{M, 2}(t) & \ldots & \hat{h}_{M, {N^\mathrm{GP} } }(t)
    \end{pmatrix}
    }}
    \cdot
    \begin{pmatrix}
      \bSS_M^1(\bx) \\
      \bSS_M^2(\bx) \\
      \vdots \\
      \bSS_M^{N{\mathrm{macro}}}(\bx)
    \end{pmatrix},
\end{equation}
where the vector of unknowns $\hat{\bh}_M(t)$ regroups the unknowns $h_e(t)$ associated with edges of conducting domain of $\Omega_c$, the unknowns $\phi_n(t)$ associated with nodes of the closure of the non-conducting domain $\overline{\Omega_c^C}$, and the currents $I_i(t)$ associated with cuts of the macroscale domain. 
The vector of shape functions $\bSS_M^i(\bx), i = 1, 2, \ldots, N^{\mathrm{macro}}$ collects $N^{\mathrm{macro}}$ macroscopic shape functions. These functions consist of Nédélec basis functions associated with the edges of the conducting domain $\Omega_c$, gradients of nodal basis functions associated with the nodes of the closure of the non-conducting domain $\overline{\Omega_c^C}$ and the global functions $\bc_i$ associated with the macroscale cuts.
Mesoscale magnetic fields are similarly approximated by: for $i = 1, 2, \cdots, N^{\mathrm{GP}}$
\begin{equation}
  \bh_c^{(i)}(\by, t) \approx 
    \underset{\hat{\bh}_c^{(i)}(t)}{\underbrace{
      \begin{pmatrix}
        \hat{h}_{c, 1}^{(i)}(t) & \hat{h}_{c, 2}^{(i)}(t) & \ldots & \hat{h}_{c, N_{\mathrm{meso}}^{(i)} }(t)
      \end{pmatrix}
      }
    }
    \cdot
    \begin{pmatrix}
      \bSS_c^{(i), 1}(\by) \\
      \bSS_c^{(i), 2}(\by) \\
      \vdots \\
      \bSS_c^{(i), N_{\mathrm{meso}}}(\by)
    \end{pmatrix},
\end{equation}
where a distinction is also made between degrees of freedom and basis functions defined in the conducting domain $\Omega_{mc}$ and in the closure of the non-conducting domain $\overline{\Omega_{mc}^C}$.
This leads to the following semi-discrete system of ordinary differential equations (ODE):
\begin{subequations}
  \begin{align}
    \begin{split}
        &\frac{d \bMM_M }{d t} \left(\hat{\bh}_M, [\hat{\bh}_c]^{N^\mathrm{GP}} \right) + \bKK_M \hat{\bh}_M = \bzero,
        \label{eq:homogenized-macro-semi-discrete}
    \end{split}
    \\
    \begin{split}
        &\frac{d \bMM_{m}^{(i)} }{d t} \left(\hat{\bh}_M + \hat{\bh}_{c}^{(i)}\right) + \bKK_m^{(i)} \hat{\bh}_{c}^{(i)} = \bzero,
        \label{eq:homogenzed-meso-semi-discrete}
    \end{split}
  \end{align}
\end{subequations}
with entries of the mass vectors $\bMM_M, \bMM_{m}^{(i)}$ and stiffness matrices $\bKK_M, \bKK_{m}^{(i)}$ defined by 
\begin{multline*}
  \left(\bMM_{M}\right)_{l} := \left(\bBB_M(\bh_M, [\bh_c]^{N^\mathrm{GP}}), \bSS_M^{l} \right)_{\Omega}, 
  \left(\bMM_{m}^{(i)}\right)_{l} := \left(\bBB(\bh_M + \bh_c^{(i)}), \bSS_c^{(i), l} \right)_{\Omega_{m}},\\
  \left(\bKK_{M}\right)_{k, l} \! := \!  \left(\rho \Curl[_x]{\bSS_M^{k}}\! , \Curl[_x]{\bSS_M^{l}} \right)_{\! \Omega_c}\! ,
  \left(\bKK_{m}^{(i)}\right)_{k, l} \! := \!  \left(\rho \Curl[_y]{\bSS_c^{(i), k}}\! , \Curl[_y]{\bSS_c^{(i), l}} \right)_{\! \Omega_{mc}}\! .
\end{multline*}

The ODEs \eqref{eq:homogenized-macro-semi-discrete}--\eqref{eq:homogenzed-meso-semi-discrete} are discretized in time using the backward Euler method at timestep $n = 1, 2, \dots, N_{\mathrm{TS}}$:
\begin{subequations}
  \begin{multline}
    \!\!\!\!\bRR_M\left(\hat{\bh}_M^{n+1}, [\hat{\bh}_c^{n+1}]^{N^\mathrm{GP}} \right) := \frac{\bMM_M(\hat{\bh}_M^{n+1}, [\hat{\bh}_c^{n+1}]^{N^\mathrm{GP}}) - \bMM_M(\hat{\bh}_M^{n}, [\hat{\bh}_c^{n}]^{N^\mathrm{GP}})}{\Delta t} \\
    + \bKK_M \hat{\bh}_M^{n+1} = \bzero,
      \label{eq:homogenized-macro-fully-discrete}
  \end{multline}
  \begin{multline}
    \!\!\!\!\bRR_{m}^{(i)}\left(\hat{\bh}_M^{n+1} + \hat{\bh}_c^{(i), n+1} \right) := \frac{\bMM_m^{(i)}(\hat{\bh}_M^{n+1} + \hat{\bh}_c^{(i), n+1} ) - \bMM_m^{(i)}(\hat{\bh}_M^{n} + [\hat{\bh}_c^{(i), n}])}{\Delta t} \\
    + \bKK_m^{(i)} \hat{\bh}_c^{(i), n+1} = \bzero.
      \label{eq:homogenized-meso-fully-discrete}
  \end{multline}
\end{subequations}
The resulting system of nonlinear equations \eqref{eq:homogenized-macro-fully-discrete}--\eqref{eq:homogenized-meso-fully-discrete} is solved using a possibly under-relaxed quasi Newton--Raphson method that reads: for $l = 1, 2, \ldots, N_\mathrm{NR}$, solve:
\begin{subequations}
  \begin{equation}
    \!\!\!\!\displaystyle \frac{\partial \bRR_M \textcolor{white}{ab} }{\partial \hat{\bh}_M^{n+1}} \underset{\Delta \hat{\bh}_M^{n+1, l+1}}{\underbrace{\bigg(\hat{\bh}_M^{n+1, l+1} - \hat{\bh}_M^{n+1, l} \bigg)}} = -\bRR_M \left(\hat{\bh}_M^{n+1, l}, [\hat{\bh}_{c}^{n+1, l} ] \right),
      \label{eq:homogenized-homog-fully-discrete_linearized_a}
  \end{equation}
  \begin{equation}
    \!\!\!\!\displaystyle \frac{\partial \bRR_{m}^{(i)} \textcolor{white}{aban} }{\partial \hat{\bh}_{c}^{(i), n+1}} \underset{\Delta \hat{\bh}_{c}^{(i), n+1, l+1}}{\underbrace{\bigg(\hat{\bh}_{c}^{(i), n+1, l+1} - \hat{\bh}_{c}^{(i), n+1, l} \bigg) } } = -\bRR_{m}^{(i)} \left(\hat{\bh}_M^{n+1, l+1} + \hat{\bh}_{c}^{(i), n+1, l} \right).
      \label{eq:homogenized-homog-fully-discrete_linearized}
  \end{equation}
\end{subequations}
with 
\begin{multline}
    \left(\frac{\partial \bRR_M \textcolor{white}{a} }{\partial \hat{\bh}_M^{n+1}} \right)_{\alpha \beta} 
    = \left( \frac{\partial \bBB_M}{\partial \hat{\bh}_M^{n+1}}\bigg(\hat{\bh}_M^{n+1, l}, [\hat{\bh}_c^{n+1, l}] \bigg) \bSS_M^{\alpha}, \bSS_M^{\beta} \right)_{\Omega} \\
    + \Delta t \left( \rho \, \Curl[_x]{\bSS_M^{\alpha}}, \Curl[_x]{\bSS_M^{\beta}} \right)_{\Omega_c}
    \label{eq:homogenized-macro-fully-discrete_linearized}
\end{multline}
and
\begin{multline}
  \left(\frac{\partial \bRR_m^{(i)} \textcolor{white}{aban} }{\partial \hat{\bh}_{c}^{(i), n+1}} \right)_{\alpha \beta}
  = \left( \frac{\partial \bBB \textcolor{white}{abana} }{\partial \hat{\bh}_{c}^{(i), n+1} } \bigg(\hat{\bh}_M^{n+1, l+1} + \hat{\bh}_{c}^{(i), n+1, l} \bigg) \bSS_{c}^{(i), \alpha}, \bSS_{c}^{(i), \beta} \right)_{\Omega{m, i}} \\
  + \Delta t \left( \rho \, \Curl[_y]{\bSS_{c}^{(i), \alpha}}, \Curl[_y]{\bSS_{c}^{(i), \beta}} \right)_{\Omega_{mc, i}}.
  \label{eq:homogenized-meso-fully-discrete_linearized}
\end{multline}
The linear systems in \eqref{eq:homogenized-homog-fully-discrete_linearized} are solved using a direct solver based on a complete LU factorization implemented in MUMPS \cite{amestoy-mumps-00}.
New solutions are then obtained as:
\begin{equation}
  \begin{aligned}
    \hat{\bh}_M^{n+1, l+1} &= \hat{\bh}_M^{n+1, l} + \omega_M \Delta \hat{\bh}_M^{n+1, l+1}, \\
    \hat{\bh}_c^{(i), n+1, l+1} &= \hat{\bh}_c^{(i), n+1, l} + \omega_m^{(i)} \Delta \hat{\bh}_c^{(i), n+1, l+1}, \quad i = 1, 2, \ldots, N^{\mathrm{GP}}.
  \end{aligned}
\end{equation}
where $\omega_M, \omega_m^{(i)} \in [0.05, 1]$ are optimal relaxation factors that minimizes the macroscale and mesoscale residuals for each non linear iteration. 
It is well known that $\bh$-conforming formulations require relaxation to convergence. In the presence of magnetic saturation, small relaxation parameters $\omega_M < 1$ and $\omega_m^{(i)} < 1$ can be used to prevent divergence or oscillatory behavior during nonlinear iterations \cite{fujiwara-nonlinear-02}. Further details on the implementation can be found in \cite{jacques-phdthesis-18, dular-thesis-23}.

Once the mesoscale problems are solved, the magnetic flux density $\bBB_M$ is upscaled by simply averaging the mesoscale magnetic flux density as 
\begin{equation}
  \bBB_M^{(i)} = \frac{1}{|\Omega_m|} \int_{\Omega_m} \bBB(\bh_M + \bh_c^{(i)} ) \text{d} \Omega_m.
    \label{eq:homogenized-b}
\end{equation}
The macroscale differential reluctivity $(\partial \bBB_M/\partial \bh_M)^{(i)}$ is upscaled using the finite difference method applied to the solution of the following magnetostatic problems: 
find $\varphi_c^{(i)} \in \Hone[_{\mathrm{per}}]{\Omega_m}$ 
such that for all $\phi_c^{(i)} \in \Hone[_{\mathrm{per}}]{\Omega_m}$ 
the following weak form holds:
\begin{equation}
  \bigg( \bBB(\bh_M + \boldsymbol{\delta}_j + \Grad[]{\varphi_c^{(i)}} ), \Grad[]{\phi_c^{(i)}} \bigg)_{\Omega_m} = 0,
  \label{eq:homogenized-meso-magsta-weak}
\end{equation}
with $\boldsymbol{\delta}_j = \delta \be_j$ for $j = 1, 2, 3$ and $i = 1, 2, \ldots, N^{\mathrm{GP}}$.

Connecting all the steps, we get the following pseudocode for the homogenized problem:
{
  \begin{center}
  \begin{algorithm}
  \small
  \caption{Pseudocode for the FE-HMM algorithm}\label{alg:FE-HMM-Macro}
  \begin{algorithmic}
  \INPUT Macroscale current sources and meshes.
  \OUTPUT Macroscale and mesoscale fields, and global quantities.
  \Procedure{Macroscale problem}{} 
      \State $t \gets 0$, init macro field $\hat{\bh}_{\mathrm{M}}^{1, 1}$ and meso fields $\hat{\bh}_c^{(i), 1, 1}$,
      \For{$(n \gets 1$ To $N_{\mathrm{TS}} )$}              \, \Comment{\emph{the time loop \textcolor{white}{abana banjye ban}  } }
      \For{$(l \gets 1$ To $N_{\mathrm{NR}}^{\mathrm{M}} )$} \, \Comment{\emph{the macroscale NR loop \textcolor{white}{abanab} } }
      \For{$(i \gets 1$ To $N^{\mathrm{GP}} )$}              \!\!\!\! \Comment{\emph{parallel solve of meso-problems} }
      \State Downscale the macroscale sources $\hat{\bh}_{\mathrm{M}}^{n+1, l}$.
      \State Call the \textbf{procedure} \textsc{Mesoscale MQS} to compute $\bBB_{\mathrm{M}}$.
      \State Call the \textbf{procedure} \textsc{Mesoscale MS} to compute $\partial \bBB_{\mathrm{M}}/ \partial \bh_{\mathrm{M}}$.
      \State Upscale $\bBB_{\mathrm{M}}^{(i)}$ and $(\partial \bBB_{\mathrm{M}}/ \partial \bh_{\mathrm{M}})^{(i)}$.
      \EndFor
      \State Assemble the matrix and the RHS, and solve the macroscale problem \eqref{eq:homogenized-homog-fully-discrete_linearized_a},
      \If{(relative residual < tol$_{\text{res}}$ or relative error on global quantities < tol$_{Q}$)} \\
        \quad \quad \quad \quad \quad \quad Exit the macroscale NR loop,
      \EndIf
      \EndFor
      \EndFor
  \EndProcedure
  \\
  \Procedure{\textsc{Mesoscale MQS} }{} 
    \State Prescribe PBC for the mesoscale field and impose the source $\hat{\bh}_M^{n+1, l+1}$,
    \For{$(j \gets 1$ To $ N_{\mathrm{NR}}^{\mathrm{m}} )$} \Comment{\emph{the mesoscale NR loop } \textcolor{white}{Ib}}
      \State Assemble the matrix and the RHS and solve the mesoscale
      \State magnetoquasistatic problem \eqref{eq:homogenized-homog-fully-discrete_linearized}.
      \State Compute $\bBB_{\mathrm{M}}$ using \eqref{eq:homogenized-b}.
    \EndFor    
  \EndProcedure
  \\
  \Procedure{\textsc{Mesoscale MS} }{} 
      \For{$(i \gets 1$ To $ 3)$} \Comment{\emph{Solve 3 meso-problems for $\partial \bBB_{\mathrm{M}}/ \partial \bh_{\mathrm{M}}$} \textcolor{white}{Ib}}
      \State Prescribe PBCs for the mesoscale field and impose the source $\hat{\bh}_M^{n+1, l+1}$,
    \For{$(j \gets 1$ To $ N_{\mathrm{NR}}^{\mathrm{m}} )$} \Comment{\emph{the mesoscale NR loop } \textcolor{white}{Ib}}
        \State Assemble the matrix and solve the mesoscale magnetostatic problems \eqref{eq:homogenized-meso-magsta-weak}.
    \EndFor    
  \EndFor
  \State Compute $\partial \bBB_{\mathrm{M}}/ \partial \bh_{\mathrm{M}}$ using finite differences.
  \EndProcedure
  \end{algorithmic}
  \end{algorithm}
  \end{center}
}
\normalsize
\clearpage
%
%
\section{Results}
\label{sec:results}
%
%
\colorbox{black!10}{$\mathbf{\Omega_{{\mathrm{c}}_{2}}}$}
This section contains the numerical tests that validate the proposed $\bh$-conforming multiscale formulations. We consider two geometries of idealized soft magnetic composites (SMCs) with a periodic cell: 
(1) a 2D SMC with conducting and magnetic disc-shaped inclusions; and 
(2) a 3D SMC with conducting and magnetic ball-shaped inclusions.

We consider the following material properties:
\begin{enumerate}
  \item The air and the insulator are non-magnetic and non-conducting with $\mu_r = 1$ and $\sigma = 0$. 
  \item The inclusions are conducting with $\sigma = 10^7$S/m, and magnetic characterised either by a linear magnetic law with a relative magnetic permeability $\mu_r \gg 1$, or by the Fr\"{o}hlich-Kennely non-linear law given by \cite{lee-materiallaw-01}: 
  \begin{equation} 
    \bb = \mu(|\bh|) \bh = \left( \mu_0 + \frac{(\mu_{\mathrm{Fe}} - \mu_0) b_s}{(\mu_{\mathrm{Fe}} - \mu_0)|\bh| +b_s} \right) \bh
  \end{equation} 
    with $\mu_{\mathrm{Fe}} = 10^2 \mu_0$ and $b_s$ = 1.5 T.
\end{enumerate}
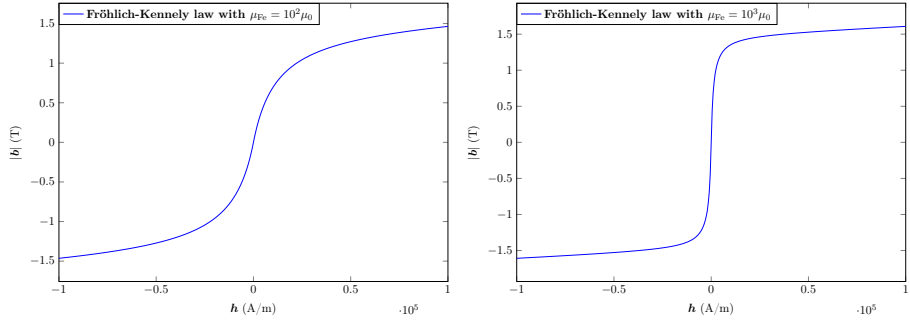
\begin{figure}
  \begin{tikzpicture}[scale=0.45]
    \begin{axis}[xlabel={$\bh$ (A/m)}, ylabel={$|\bb|$ (T)}, xmin=-100000.0, xmax=100000.0, xtick={-100000, -50000, 0, 50000, 100000}, width=1.0\columnwidth, height=0.75\columnwidth, legend style={at={(0,1)},anchor=north west}]
    \addplot[domain = -100000:100000, samples = 200, smooth, thick, blue] 
      {( 4*pi*1.0e-7 + (4*pi*1.0e-7*(100-1) * 1.5)/(4*pi*1.0e-7*(100-1)*abs(x) + 1.5 ) ) * x};
      \addlegendentry{\textbf{Fr\"{o}hlich-Kennely law with $\mu_{\mathrm{Fe}} = 10^2 \mu_0$}}
      \end{axis}
  \end{tikzpicture}
  \begin{tikzpicture}[scale=0.45]
    \begin{axis}[xlabel={$\bh$ (A/m)}, ylabel={$|\bb|$ (T)}, xmin=-100000.0, xmax=100000.0, xtick={-100000, -50000, 0, 50000, 100000}, width=1.0\columnwidth, height=0.75\columnwidth, legend style={at={(0,1)},anchor=north west}]
    \addplot[domain = -100000:100000, samples = 200, smooth, thick, blue] 
      {( 4*pi*1.0e-7 + (4*pi*1.0e-7*(1000-1) * 1.5)/(4*pi*1.0e-7*(1000-1)*abs(x) + 1.5 ) ) * x};
      \addlegendentry{\textbf{Fr\"{o}hlich-Kennely law with $\mu_{\mathrm{Fe}} = 10^3 \mu_0$}}
      \end{axis}
  \end{tikzpicture}
  \caption{\footnotesize
  Two examples of the Fr\"{o}hlich-Kennely non linear laws with $b_s = 1.5$ T. 
  \textbf{Left} : BH curve used for 3D problems with the permeability at the origin $\mu_{\mathrm{Fe}} = 10^2 \mu_0$. 
  \textbf{Right} : BH curve used for 2D problems with the permeability at the origin $\mu_{\mathrm{Fe}} = 10^3 \mu_0$.}
    \label{HMM_2D_Geo_Meshes}
\end{figure}
\normalsize
The problems can be fed by a sinusoidal voltage or current sources as detailed in Section~\ref{sec:finescale_formulations}. In this paper, we only present results with the sinusoidal current sources $I_{s}(t) = I_{s0} \sin{(2 \pi f t)}$ where $I_{s0}$ will be determined for each case. Simulations are carried out using Gmsh~\cite{geuzaine-gmsh-09} for the generation of geometries and meshes, as well as for postprocessing, and GetDP~\cite{dular-getdp-98} for numerical resolutions. 
Parallelization was carried out using the supercomputer NIC5 of the CECI cluster~\cite{nic5-ceci-1} and on the supercomputer of the GRICAD of Université Grenoble Alpes \cite{gricad-uga-1}.

For our validation, we focus on the accuracy of the models with eddy currents and nonlinear material laws in 2D and 3D, and on the performance of the parallel algorithms in 3D.
For post-processing, we define the homogenized eddy-current losses, the magnetic power and the total power in the homogenized core, respectively as:
\begin{multline}
  \mathrm{P}_{\Omega_{\mathrm{core}}}^{\mathrm{HMM, mag}} = \int_{\Omega_{\mathrm{core}}} \mathrm{P}_{\mathrm{m}}^{\mathrm{mag}} \text{d}\Omega, \quad 
  \mathrm{P}_{\Omega_{\mathrm{core}}}^{\mathrm{HMM, eddy}} = \int_{\Omega_{\mathrm{core}}} \mathrm{P}_{\mathrm{m}}^{\mathrm{eddy}} \text{d}\Omega \\
  \quad \text{ and } \quad    
  \mathrm{P}_{\Omega_{\mathrm{core}}}^{\mathrm{HMM, tot}} = \int_{\Omega_{\mathrm{core}}} (\mathrm{P}_{\mathrm{m}}^{\mathrm{mag}} + \mathrm{P}_{\mathrm{m}}^{\mathrm{eddy}} ) \text{d}\Omega,
\end{multline}
where $\mathrm{P}_{\mathrm{m}}^{\mathrm{mag}}$ and $\mathrm{P}_{\mathrm{m}}^{\mathrm{eddy}}$ are average power defined in \eqref{eq:average_meso_quantities}.
Throughout this section, we define also the relative error time function for a quantity Q (where Q represents either Joule losses , the magnetic power or the voltage) as follows:
\begin{equation} 
  \varepsilon_{Q}(t) 
  = \frac{|Q^{\mathrm{ref}}(t) - Q^{\mathrm{HMM}}(t) |}{\|Q^{\mathrm{ref}}\|_{L^{\infty}(0, T)}}
  = \frac{|Q^{\mathrm{ref}}(t) - Q^{\mathrm{HMM}}(t) |}{\underset{t \in ]0, T[}{\text{max}} \quad |Q^{\mathrm{ref}} (t)|}.
  \label{eq:homog_quantities}
\end{equation}
%
%
\subsection{Two-dimensional SMCs: accuracy}
%
%
Reference results are obtained using a brute force approach with the FE problem solved on the fine mesh with 184685 elements depicted in Figure~\ref{HMM_2D_Geo_Meshes} (left).
This mesh is generated from a geometry of an idealized 2D SMC with $8 \times 8$ disks, with the period $\Omega_m$ of width 100 $\mu$m and a conducting inclusion of radius 40 $\mu$m. Homogenized results are obtained using the macroscale mesh in Figure~\ref{HMM_2D_Geo_Meshes} (middle) with 64 macro-elements in the homogenized domain, and the mesoscale mesh in Figure~\ref{HMM_2D_Geo_Meshes} (right) with 7099 elements for the cell problems. 
\begin{figure}
  \begin{center}
  \includegraphics[width=.25\textwidth]{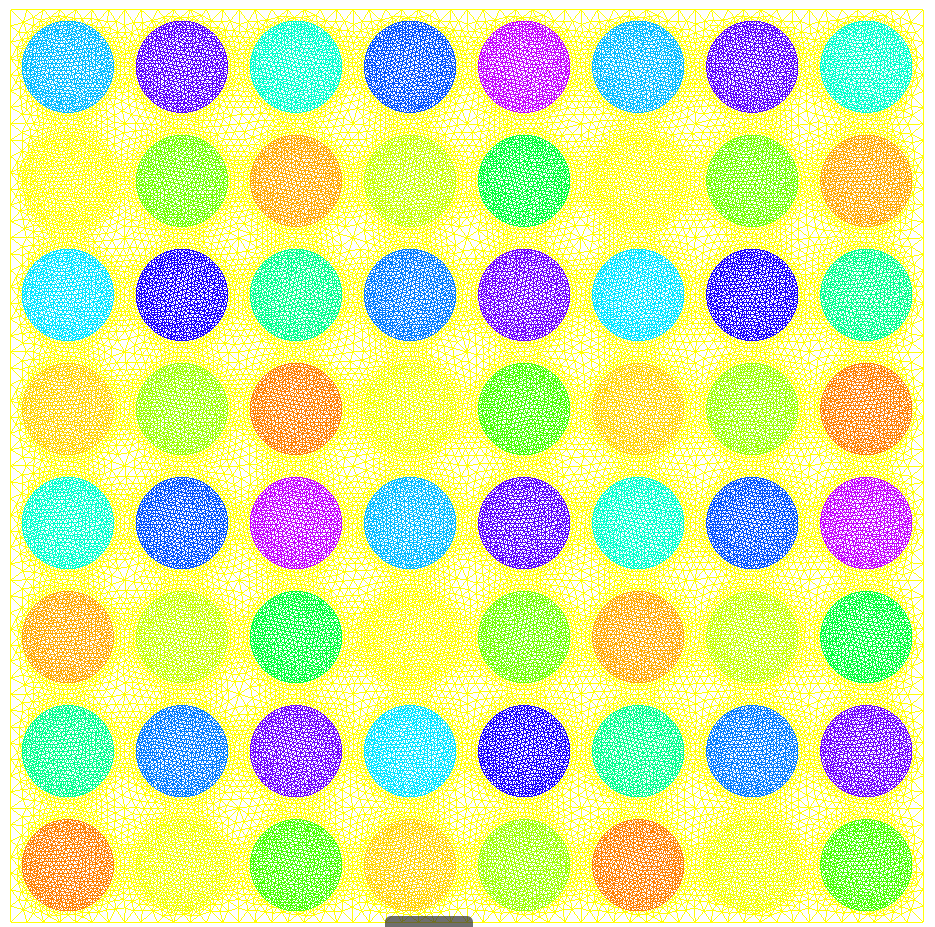}
  \hspace{0mm}
  \includegraphics[width=.25\textwidth]{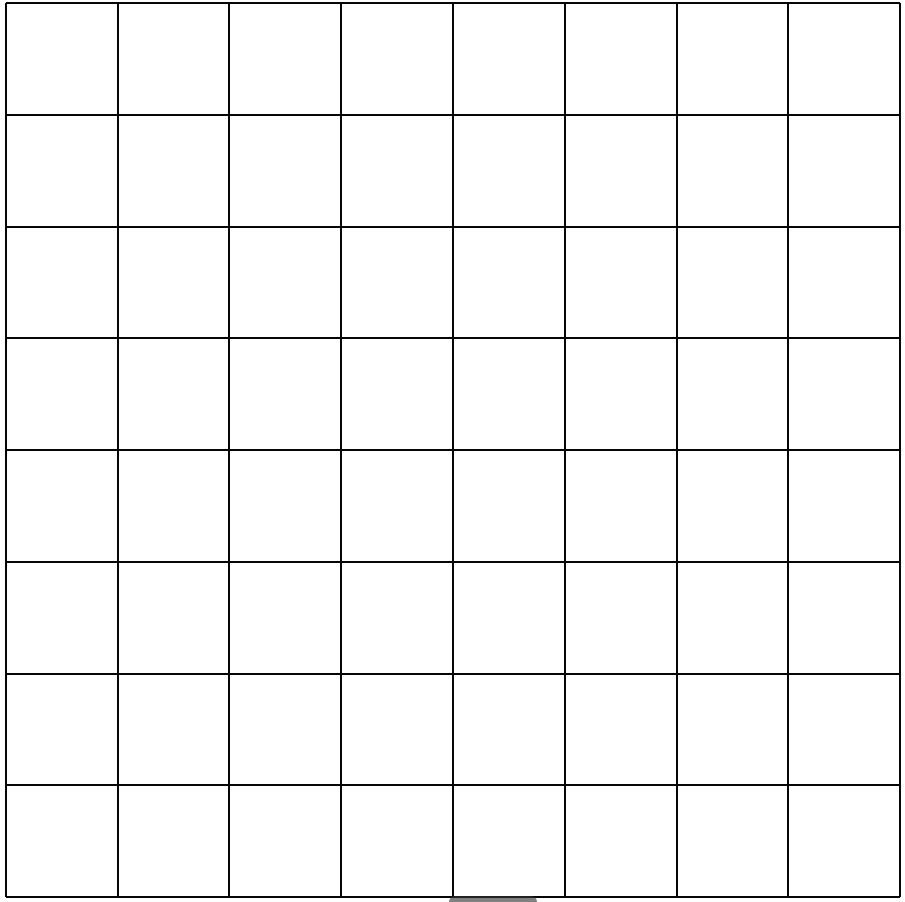}
  \hspace{0mm}
  \includegraphics[width=.25\textwidth]{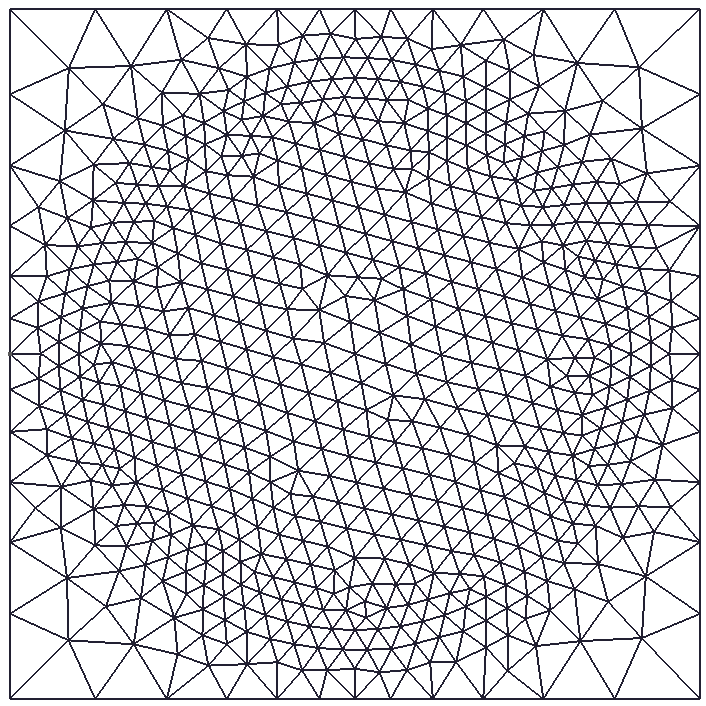}
  \end{center}
  \caption{\footnotesize
  Geometry and meshes used for the 2D problem. Left: the fullscale core. Middle: The mesh of the homogenized core. Right: The mesh of a mesoscale cell.}
  \label{HMM_2D_Geo_Meshes}
\end{figure} 
\normalsize

\begin{figure}
  \begin{center}
    \includegraphics[width=0.3\textwidth]{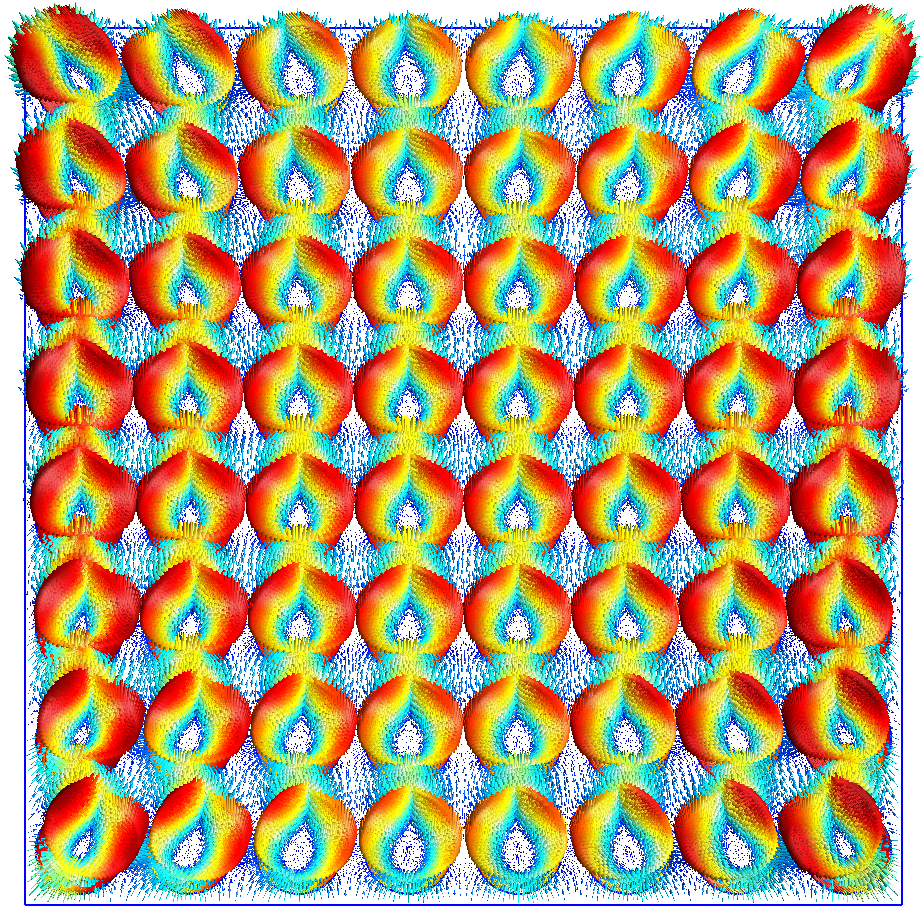}
    \includegraphics[width=0.3\textwidth]{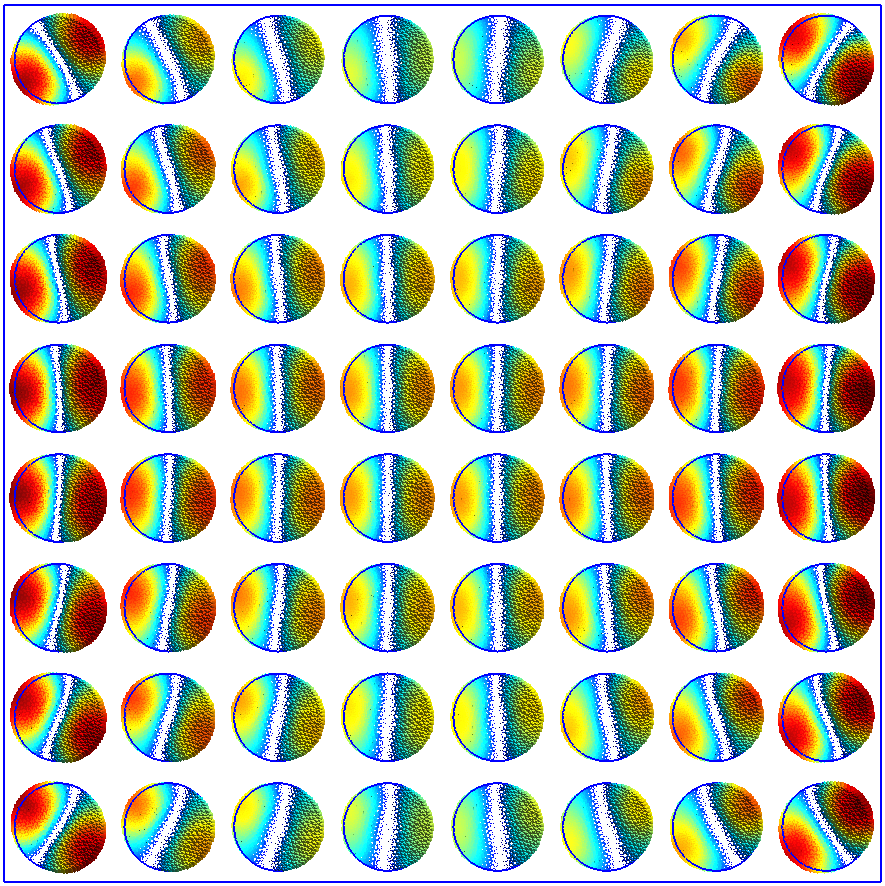}
    \end{center}
    \caption{\footnotesize
    Typical map of the reference magnetic induction $\bb^{\varepsilon}$ (left) and eddy currents $\bj^{\varepsilon}$ (right) in the core for a problem with $\mu_r^{\text{core}} = 10^3$ and the frequency $f = 100$ kHz.}
    \label{HMM_2D_b_map}
\end{figure}

Figure \ref{HMM_2D_b_map} shows typical maps of the reference magnetic induction $\bb^{\varepsilon}$ and eddy currents $\bj^{\varepsilon}$ in the core.

\normalsize 

Figures~\ref{fig:HMM_2D_lin_localquantities_1} and \ref{fig:HMM_2D_nonlin_localquantities_1} illustrate the time evolution of several global quantities (eddy-current losses, the magnetic power, the magnetic flux density $\bBB_M$ in the homogenized core, and the induced voltage), together with the corresponding accuracy of these quantities as functions of time. In both cases, the associated problems highlight the intrinsic complexity of homogenizing magnetic problems involving eddy currents and nonlinear material behavior.  

Figure \ref{fig:HMM_2D_lin_localquantities_1} is obtained for a linear non-magnetic problem ($\mu_r = 1$) at $f = 100$ MHz. The excitation is a sinusoidal current $I_{s}(t) = I_{s0} \sin{(2 \pi f t)}$ with amplitude $I_{s0} = 2000$A, and each period is discretized using 120 time steps. Figure \ref{fig:HMM_2D_nonlin_localquantities_1} is obtained solving a nonlinear problem with the Fr\"{o}hlich-Kennely nonlinear law, with $\mu_{\mathrm{Fe}} = 10^3 \mu_0$ at a frequency $f = 100$ kHz. The excitation is a sinuoidal source current $I_{s}(t) = I_{s0} \sin{(2 \pi f t)}$ with amplitude $I_{s0} = 1000$A, and each period is discretized using 160 time steps.

Several conclusions can be drawn from these figures. First of all, as it has already been mentionned in Section \ref{sec:dimensional_asymptotic_analysis_2}, the homogenized magnetic induction $\bBB_M$ obtained from the magnetostatic and magnetodynamic formulations exhibits very similar values when the ratio of eddy-current losses to the total power within the cell is small (see Figure \ref{fig:HMM_2D_nonlin_localquantities_1} Middle-left). In contrast, noticeable discrepancies appear when this ratio becomes non-negligible (see Figure \ref{fig:HMM_2D_lin_localquantities_1} Middle-left). Furthermore, the eddy-current losses and the magnetic power in the core, as well as the induced voltage are accurately captured within the HMM framework. In particular, the relative errors in the eddy-current losses remain below 0.2 \% for the linear case and below 5 \% for the nonlinear case. For all test cases considered in this paper, the absolute error profiles exhibit the same temporal behavior as the corresponding relative errors, differing only by a scaling factor. For this reason, absolute error curves are shown only for the two-dimensional linear case. The absolute errors for the other configurations can be readily obtained by multiplying the relative errors by the maximum value of the corresponding quantity of interest which for the linear problem are approximately $2 \times 10^{5}$ W for the eddy-current losses and $4 \times 10^{6}$ V for the voltage. Similarly, the relative error in the voltage stays below 5 \% for the linear problem and below 0.003 \% for the nonlinear case. 

The stopping criteria were set to a relative residual of $10^{-8}$ for the reference problem and $10^{-4}$ for the macroscale problem. The latter threshold proved difficult to attain, as the residual of the macroscale problem converges very slowly. More effective convergence indicators were found to be the relative increment of the unknown field, as well as the relative residuals of the homogenized Joule losses and the magnetic power in the homogenized domain, as defined in \eqref{eq:homog_quantities}. These alternative criteria were therefore set to $10^{-7}$. Overall, for all multiscale simulations, we observed a significantly slower convergence of the macroscale residual compared with what is typically expected for a single-scale nonlinear problem solved using the $\bh$-conforming formulation.

The stopping criteria were set to a relative residual of $10^{-8}$ for the reference and the mesoscale problems, and $10^{-4}$ for the macroscale problem. The relative residual $\br_{\mathrm{rel}}^{n+1, l} = \br_{\mathrm{abs}}^{n+1, l}/\br_{\mathrm{abs}}^{n+1, 1}$ for the time step $n+1$ and the nonlinear iteration $l$ is defined with respect to the initial absolute residual $\br_{\mathrm{abs}}^{n+1, 1} = \| \bb - \bAA \bx^{n+1, 1} \|$. The latter threshold proved difficult to attain, as the residual of the macroscale problem converges very slowly. More effective convergence indicators were found to be the relative increment of the unknown field, as well as the relative residuals of the homogenized Joule losses and the magnetic power in the homogenized domain. These alternative criteria were therefore set to $10^{-7}$. Overall, for all multiscale simulations, we observed a significantly slower convergence of the macroscale residual compared with what is typically expected for a single-scale nonlinear problem solved using the $\bh$-conforming formulation.

\begin{figure}[H]
  \begin{tikzpicture}[scale=0.5]
    \begin{axis}[xlabel={Time (s)}, ylabel={Joule Losses (W)}, xmin=0.0, xmax=2e-8, xtick={0, 0.5e-8, 1e-8, 1.5e-8, 2e-8}, ymin=0.0, 
      width=0.98\columnwidth, height=0.5\columnwidth, legend style={at={(0,1)},anchor=north west}]
      \addplot [color=blue, mark=none, line width=1.0mm, mark size=2.0, mark options=solid] table {./data/2D_smc_ref_lin_2e3A_f1e8_mur1_NSteps120_JouleLosses_Core_Ref.csv};
      \addlegendentry{\textbf{Ref}}
      \addplot [color=red, mark=*, mark size=3.0, mark options=dashdotted] table {./data/2D_smc_hmm_lin_2e3A_f1e8_mur1_NSteps120_JL_Homog.csv};
      \addlegendentry{\textbf{HMM}}
    \end{axis}
  \end{tikzpicture}
  \begin{tikzpicture}[scale=0.5]
    \begin{axis}[xlabel={Time (s)}, ylabel={Magnetic power (W)}, xmin=0.0, xmax=2e-8, xtick={0, 0.5e-8, 1e-8, 1.5e-8, 2e-8}, 
      width=0.98\columnwidth, height=0.5\columnwidth, legend style={at={(0,1)},anchor=north west}]
      \addplot [color=blue, mark=none, line width=1.0mm, mark size=2.0, mark options=solid] table {./data/2D_smc_ref_lin_2e3A_f1e8_mur1_NSteps120_MagPower_Ref.csv};
      \addlegendentry{\textbf{Ref}}
      \addplot [color=red, mark=*, mark size=3.0, mark options=dashdotted] table {./data/2D_smc_hmm_lin_2e3A_f1e8_mur1_NSteps120_MP_Integral.csv};
      \addlegendentry{\textbf{HMM}}
    \end{axis}
  \end{tikzpicture}

  \begin{tikzpicture}[scale=0.5]
    \begin{axis}[xlabel={Time (s)}, ylabel={Upscaled magnetic induction ${\bBB_M}_y$ (T)}, xmin=0.0, xmax=2e-8, xtick={0, 0.5e-8, 1e-8, 1.5e-8, 2e-8}, 
      width=0.95\columnwidth, height=0.5\columnwidth, legend style={at={(0,1)},anchor=north west}]
      \addplot [color=blue, mark=none, line width=1.0mm, mark size=3.0, mark options=solid] table[x index=0,y index=2,col sep=space] {./data/2D_smc_hmm_lin_2e3A_f1e8_mur1_NSteps120_b_mean_GP7952.csv};
      \addlegendentry{\textbf{Magnetodynamics}}
      \addplot [color=red, mark=*, line width=1.0mm, mark size=3.0, mark options=solid] table[x index=0,y index=2,col sep=space] {./data/2D_smc_hmm_lin_2e3A_f1e8_mur1_NSteps120_b_mean_MagSta_GP7952.csv};
      \addlegendentry{\textbf{Magnetostatics}}
    \end{axis}
  \end{tikzpicture}
  \begin{tikzpicture}[scale=0.5]
    \begin{axis}[xlabel={Time (s)}, ylabel={Voltage U (V)}, axis y line*=left, xmin=0.0, xmax=2e-8, xtick={0, 0.5e-8, 1e-8, 1.5e-8, 2e-8}, width=0.9\columnwidth, height=0.5\columnwidth, legend style={at={(0,1)},anchor=north west}]
      \addplot [color=blue, mark=*, line width=1.0mm, mark size=2.0, mark options=solid] table {./data/2D_smc_ref_lin_2e3A_f1e8_mur1_NSteps120_U.csv};
      \addlegendentry{\textbf{Ref - Voltage}}     
      \addplot [color=red, mark=*, mark size=2.0, mark options=dashdotted] table {./data/2D_smc_hmm_lin_2e3A_f1e8_mur1_NSteps120_U_Macro.csv};
      \addlegendentry{\textbf{HMM - Voltage}}
    \end{axis}
    \begin{axis}[xlabel={Time (s)}, ylabel={Current I (A)}, axis y line*=right, xmin=0.0, xmax=2e-8, xtick={0, 0.5e-8, 1e-8, 1.5e-8, 2e-8}, ymin=-4000, ymax=4000, width=0.9\columnwidth, height=0.5\columnwidth, legend style={at={(1,1)},anchor=north east}]
      \addplot [color=black, mark=*, line width=1.0mm, mark size=2.0, mark options=solid] table [x expr=\thisrowno{0}*1, y expr=\thisrowno{1}*-1, col sep=space] {./data/2D_smc_ref_lin_2e3A_f1e8_mur1_NSteps120_I.csv};
      \addlegendentry{\textbf{Ref - Current}}     
      \addplot [color=green, mark=*, mark size=2.0, mark options=dashdotted] table [x expr=\thisrowno{0}*1, y expr=\thisrowno{1}*-1, col sep=space] {./data/2D_smc_hmm_lin_2e3A_f1e8_mur1_NSteps120_I_Macro.csv};
      \addlegendentry{\textbf{HMM - Current}}
    \end{axis}
  \end{tikzpicture}
  \begin{tikzpicture}[scale=0.5]
    \begin{axis}[xlabel={Time (s)}, ylabel={Absolute error $\varepsilon_{\mathrm{JL}}$ (W)}, axis y line*=left, xmin=0.0, xmax=2e-8, xtick={0, 0.5e-8, 1e-8, 1.5e-8, 2e-8}, width=0.9\columnwidth, height=0.5\columnwidth, legend style={at={(0,1)},anchor=north west}]
      \pgfplotstableread{./data/2D_smc_ref_lin_2e3A_f1e8_mur1_NSteps120_JouleLosses_Core_Ref.csv}{\datatable}
      \pgfplotstablecreatecol[copy column from table={./data/2D_smc_hmm_lin_2e3A_f1e8_mur1_NSteps120_JL_Homog.csv}{1}] {dataA} {\datatable}
      \addplot [color=blue, mark=none, line width=3.0mm, mark size=3.0, mark options=solid] table [x expr=\thisrowno{0}*1e0, y expr=abs( (\thisrowno{1} - \thisrow{dataA}))] {\datatable};
      \addlegendentry{\textbf{Absolute error}}     
    \end{axis}
    \begin{axis}[xlabel={Time (s)}, ylabel={Relative error $\varepsilon_{\mathrm{JL}}$ (\%)}, axis y line*=right, xmin=0.0, xmax=2e-8, xtick={0, 0.5e-8, 1e-8, 1.5e-8, 2e-8}, width=0.9\columnwidth, height=0.5\columnwidth, legend style={at={(1,1)},anchor=north east}]
      \pgfplotstableread{./data/2D_smc_ref_lin_2e3A_f1e8_mur1_NSteps120_JouleLosses_Core_Ref.csv}{\datatable}
      \pgfplotstablecreatecol[copy column from table={./data/2D_smc_hmm_lin_2e3A_f1e8_mur1_NSteps120_JL_Homog.csv}{1}] {dataA} {\datatable}
      \addplot [color=red, mark=*, line width=1.0mm, mark size=2.0, mark options=solid] table [x expr=\thisrowno{0}*1e0, y expr=abs( (\thisrowno{1}/99458808 - \thisrow{dataA}/99458808))] {\datatable};
      \addlegendentry{\textbf{Relative error}}     
    \end{axis}
  \end{tikzpicture}
  \begin{tikzpicture}[scale=0.5]
    \begin{axis}[xlabel={Time (s)}, ylabel={Absolute error $\varepsilon_{v}$}, axis y line*=left, xmin=0.0, xmax=2e-8, xtick={0, 0.5e-8, 1e-8, 1.5e-8, 2e-8}, width=0.9\columnwidth, height=0.5\columnwidth, legend style={at={(0,1)},anchor=north west}]
      \pgfplotstableread{./data/2D_smc_ref_lin_2e3A_f1e8_mur1_NSteps120_U.csv}{\datatable}
      \pgfplotstablecreatecol[copy column from table={./data/2D_smc_hmm_lin_2e3A_f1e8_mur1_NSteps120_U_Macro.csv}{1}] {dataA} {\datatable}
      \addplot [color=blue, mark=nano, line width=3.0mm, mark size=3.0, mark options=solid] table [x expr=\thisrowno{0}*1e0, y expr=abs( (\thisrowno{1} - \thisrow{dataA}))] {\datatable};
      \addlegendentry{\textbf{Absolute error}}     
    \end{axis}
    \begin{axis}[xlabel={Time (s)}, ylabel={Relative error $\varepsilon_{v}$ (\%)}, axis y line*=right, xmin=0.0, xmax=2e-8, xtick={0, 0.5e-8, 1e-8, 1.5e-8, 2e-8}, width=0.9\columnwidth, height=0.5\columnwidth, legend style={at={(0,1)},anchor=north west}]
      \pgfplotstableread{./data/2D_smc_ref_lin_2e3A_f1e8_mur1_NSteps120_U.csv}{\datatable}
      \pgfplotstablecreatecol[copy column from table={./data/2D_smc_hmm_lin_2e3A_f1e8_mur1_NSteps120_U_Macro.csv}{1}] {dataA} {\datatable}
      \addplot [color=red, mark=*, line width=1.0mm, mark size=2.0, mark options=solid] table [x expr=\thisrowno{0}*1e0, y expr=abs( (\thisrowno{1}/38288667 - \thisrow{dataA}/38288667))] {\datatable};
      \addlegendentry{\textbf{Relative error}}     
    \end{axis}
  \end{tikzpicture}
  \caption{\footnotesize
    Global quantities for the linear non-magnetic problem with ($\mu_r = 1$) at $f = 100$ MHz. 
    \textbf{Top left}: Joule losses in the core obtained solving the reference problem (blue curve) and the homogenized problem (red curve). 
    \textbf{Top right}: Magnetic power in the core obtained solving the reference problem (blue curve) and the homogenized problem (red curve).
    \textbf{Middle left}: The $\by$-component of the homogenized magnetic field ${\bBB_M}_y$ at one Gauss point of the macro-element obtained solving magnetodynamic problem (blue curve) and magnetostatic problem (red curve). 
    \textbf{Middle right}: Homogenized voltage (red curve) and current (red curve).
    \textbf{Bottom left}: Relative error $\varepsilon_{\mathrm{JL}}$ on Joule losses. 
    \textbf{Bottom right}: Relative error $\varepsilon_{v}$ on voltage. 
  }
  \label{fig:HMM_2D_lin_localquantities_1}
\end{figure}
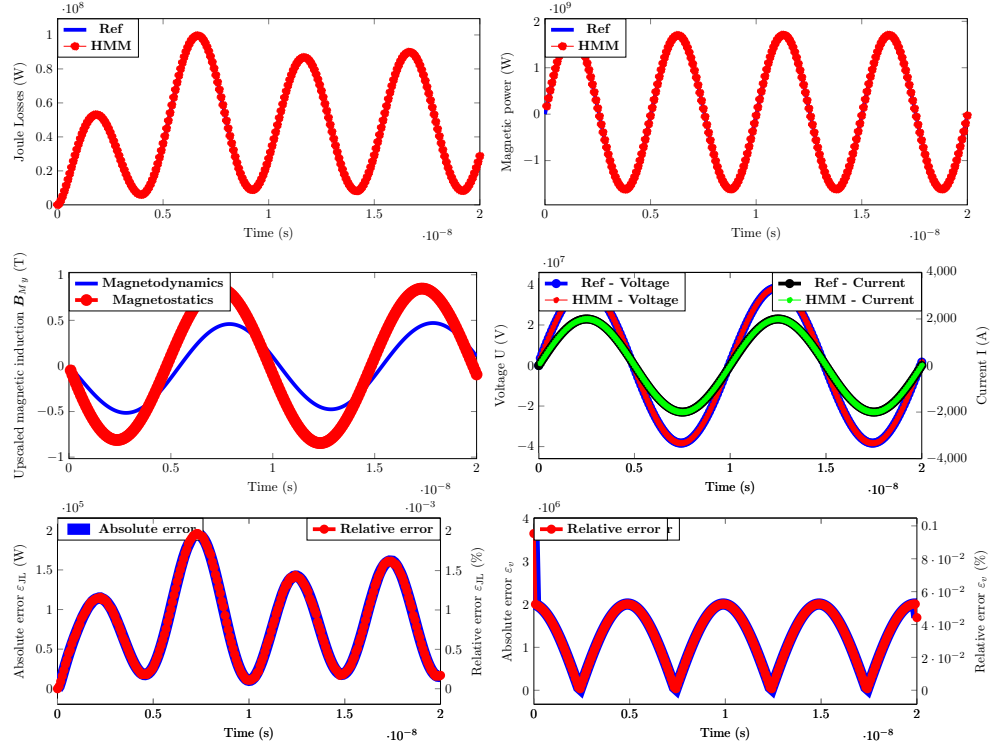
\begin{figure}[H]
  \begin{tikzpicture}[scale=0.5]
    \begin{axis}[xlabel={Time (s)}, ylabel={Joule Losses (W)}, xmin=0.0, xmax=2.e-5, xtick={0, 0.5e-5, 1.0e-5, 1.5e-5, 2.0e-5}, ymin=0.0,
      width=0.98\columnwidth, height=0.5\columnwidth, legend style={at={(0,1)},anchor=north west}]
      \addplot [color=blue, mark=none, line width=1.0mm, mark size=2.0, mark options=solid] table {./data/2D_smc_ref_nl_fk_1e3A_f1e5_mur1000_NSteps160_JL_Core_Ref.csv};
      \addlegendentry{\textbf{Ref}}
      \addplot [color=red, mark=*, mark size=2.0, mark options=dashdotted] table {./data/2D_smc_hmm_nl_fk_1e3A_f1e5_mur1000_NSteps320_JL_Homog.csv};
      \addlegendentry{\textbf{HMM}}
    \end{axis}
  \end{tikzpicture}
  \begin{tikzpicture}[scale=0.5]
    \begin{axis}[xlabel={Time (s)}, ylabel={Magnetic power (W)}, xmin=0.0, xmax=2e-5, xtick={0, 0.5e-5, 1e-5, 1.5e-5, 2e-5}, 
      width=0.98\columnwidth, height=0.5\columnwidth, legend style={at={(0,1)},anchor=north west}]
      \addplot [color=blue, mark=none, line width=1.0mm, mark size=2.0, mark options=solid] table {./data/2D_smc_ref_nl_fk_1e3A_f1e5_mur1000_NSteps160_MagPower_Ref.csv};
      \addlegendentry{\textbf{Ref}}
      \addplot [color=red, mark=*, mark size=2.0, mark options=dashdotted] table {./data/2D_smc_hmm_nl_fk_1e3A_f1e5_mur1000_NSteps320_MP_Integral.csv};
      \addlegendentry{\textbf{HMM}}
    \end{axis}
  \end{tikzpicture}
  \begin{tikzpicture}[scale=0.5]
    \begin{axis}[xlabel={Time (s)}, ylabel={Upscaled magnetic induction ${\bBB_M}_y$ (T)}, xmin=0.0, xmax=2e-5, xtick={0, 0.5e-5, 1e-5, 1.5e-5, 2e-5}, 
      width=0.95\columnwidth, height=0.5\columnwidth, legend style={at={(0,1)},anchor=north west}]
      \addplot [color=blue, mark=*, line width=1.0mm, mark size=2.0, mark options=solid] table[x index=0,y index=2,col sep=space] {./data/2D_smc_hmm_nl_fk_1e3A_f1e5_mur1000_NSteps320_b_mean_GP7346.csv};
      \addlegendentry{\textbf{Magnetodynamics}}
      \addplot [color=red, mark=*, line width=1.0mm, mark size=2.0, mark options=solid] table[x index=0,y index=2,col sep=space] {./data/2D_smc_hmm_nl_fk_1e3A_f1e5_mur1000_NSteps320_b_MS_GP7346.csv};
      \addlegendentry{\textbf{Magnetostatics}}
    \end{axis}
  \end{tikzpicture}
  \begin{tikzpicture}[scale=0.5]
    \begin{axis}[xlabel={Time (s)}, ylabel={Voltage U (V)}, axis y line*=left, xmin=0.0, xmax=2e-5, xtick={0, 0.5e-5, 1e-5, 1.5e-5, 2e-5}, width=0.90\columnwidth, height=0.5\columnwidth, legend style={at={(0,1)},anchor=north west}]
      \addplot [color=blue, mark=*, line width=1.0mm, mark size=2.0, mark options=solid] table {./data/2D_smc_ref_nl_fk_1e3A_f1e5_mur1000_NSteps160_U.csv};
      \addlegendentry{\textbf{Ref - Voltage}}     
      \addplot [color=red, mark=*, mark size=2.0, mark options=dashdotted] table {./data/2D_smc_hmm_nl_fk_1e3A_f1e5_mur1000_NSteps320_U_Macro.csv};
      \addlegendentry{\textbf{HMM - Voltage}}
    \end{axis}
    \begin{axis}[xlabel={Time (s)}, ylabel={Current I (A)}, axis y line*=right, xmin=0.0, xmax=2e-5, xtick={0, 0.5e-5, 1e-5, 1.5e-5, 2e-5}, ymin=-2000, ymax=2000, width=0.90\columnwidth, height=0.5\columnwidth, legend style={at={(1,1)},anchor=north east}]
      \addplot [color=black, mark=*, line width=1.0mm, mark size=2.0, mark options=solid] table [x expr=\thisrowno{0}*1, y expr=\thisrowno{1}*-1, col sep=space] {./data/2D_smc_ref_nl_fk_1e3A_f1e5_mur1000_NSteps160_I.csv};
      \addlegendentry{\textbf{Ref - Current}}     
      \addplot [color=green, mark=*, mark size=2.0, mark options=dashdotted] table [x expr=\thisrowno{0}*1, y expr=\thisrowno{1}*-1, col sep=space] {./data/2D_smc_hmm_nl_fk_1e3A_f1e5_mur1000_NSteps320_I_Macro.csv};
      \addlegendentry{\textbf{HMM - Current}}
    \end{axis}
  \end{tikzpicture}
  \begin{tikzpicture}[scale=0.5]
    \begin{axis}[xlabel={Time (s)}, ylabel={Relative error $\varepsilon_{\mathrm{JL}}$}, xmin=0.0, xmax=2.e-5, xtick={0, 0.5e-5, 1.0e-5, 1.5e-5, 2.0e-5}, width=0.97\columnwidth, height=0.5\columnwidth, legend style={at={(0,1)},anchor=north west}]
      \pgfplotstableread{./data/2D_smc_ref_nl_fk_1e3A_f1e5_mur1000_NSteps160_JL_Core_Ref.csv}{\datatable}
      \pgfplotstablecreatecol[copy column from table={./data/2D_smc_hmm_nl_fk_1e3A_f1e5_mur1000_NSteps320_JL_Homog_Halved.csv}{1}] {dataA} {\datatable}
      \addplot [color=red, mark=*, line width=1.0mm, mark size=2.0, mark options=solid] table [x expr=\thisrowno{0}*1e0, y expr=abs( (\thisrowno{1}/542 - \thisrow{dataA}/542))] {\datatable};
    \end{axis}
  \end{tikzpicture}
  \begin{tikzpicture}[scale=0.5]
    \begin{axis}[xlabel={Time (s)}, ylabel={Relative error $\varepsilon_{v}$}, xmin=0.0, xmax=2.e-5, xtick={0, 0.5e-5, 1.0e-5, 1.5e-5, 2.0e-5}, width=0.97\columnwidth, height=0.5\columnwidth, legend style={at={(0,1)},anchor=north west}]
      \pgfplotstableread{./data/2D_smc_ref_nl_fk_1e3A_f1e5_mur1000_NSteps160_U.csv}{\datatable}
      \pgfplotstablecreatecol[copy column from table={./data/2D_smc_hmm_nl_fk_1e3A_f1e5_mur1000_NSteps320_U_Macro_Halved.csv}{1}] {dataA} {\datatable}
      \addplot [color=red, mark=*, line width=1.0mm, mark size=2.0, mark options=solid] table [x expr=\thisrowno{0}*1e0, y expr=abs( (\thisrowno{1}/19079977 - \thisrow{dataA}/19079977))] {\datatable};
    \end{axis}
  \end{tikzpicture}
  \caption{\footnotesize
    Global quantities for the nonlinear problem at $f = 100$ kHz. 
    \textbf{Top left}: Joule losses obtained solving the reference problem (blue curve) and the homogenized problem (red curve). 
    \textbf{Top right}: Magnetic power obtained solving the reference problem (blue curve) and the homogenized problem (red curve).
    \textbf{Middle left}: The y-component of the homogenized magnetic field ${\bBB_M}_y$ at a Gauss point obtained solving magnetodynamic problem (blue curve) and magnetostatic problem (red curve). 
    \textbf{Middle right}: Homogenized voltage (red curve) and current (red curve).
    \textbf{Bottom left}: Relative error $\varepsilon_{\mathrm{JL}}$ on Joule losses. 
    \textbf{Bottom right}: Relative error $\varepsilon_{v}$ on voltage. 
  }
  \label{fig:HMM_2D_nonlin_localquantities_1}
\end{figure}
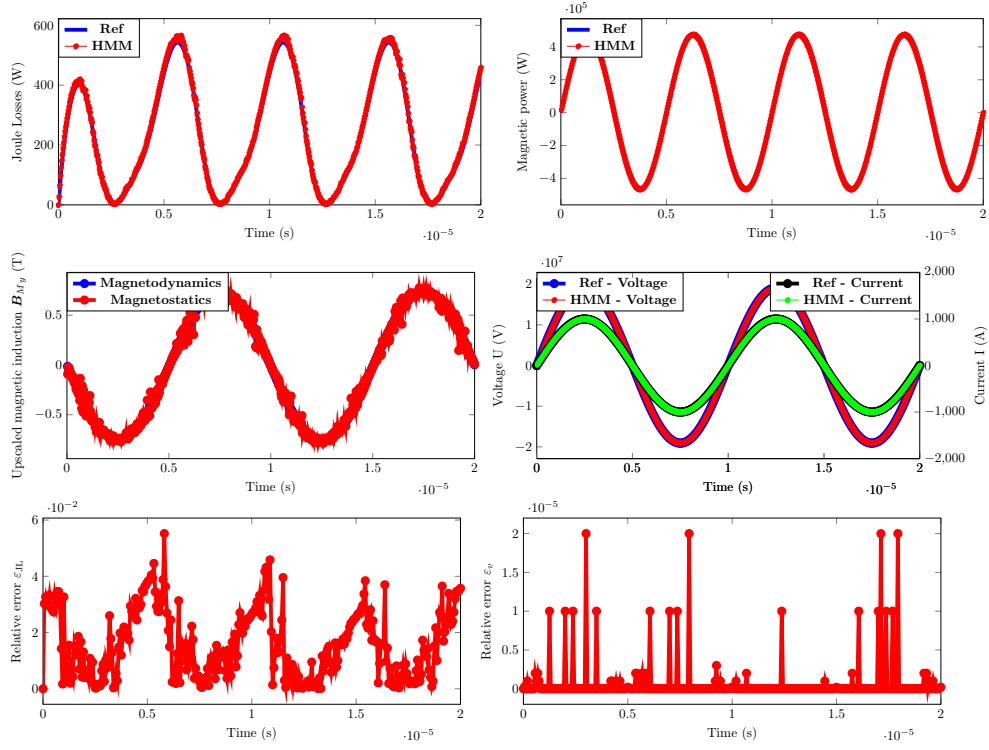
\normalsize
%
%
%
%
\subsection{Three-dimensional SMCs: accuracy}
\begin{figure}[H]
    \begin{center}
      \includegraphics[scale=0.18]{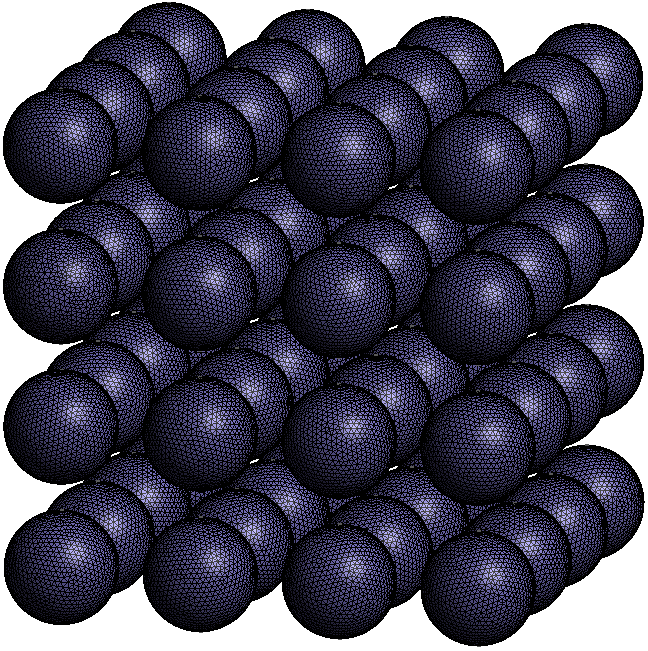} 
      \hspace{0mm}
      \includegraphics[scale=0.19]{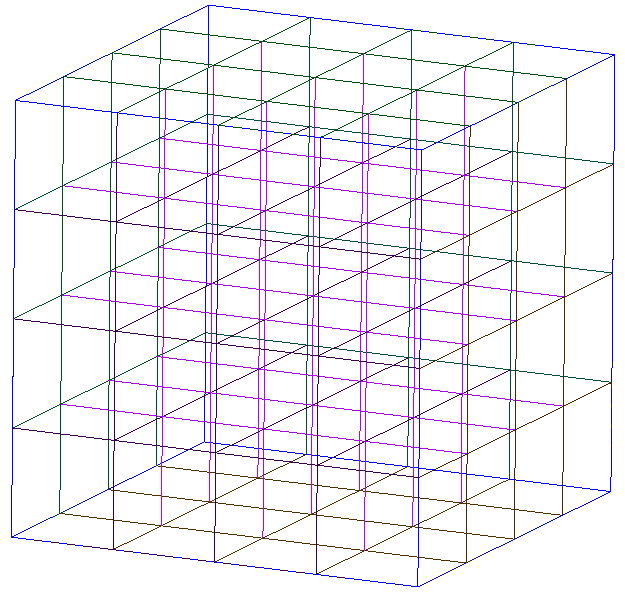}
      \hspace{0mm}
      \includegraphics[scale=0.11]{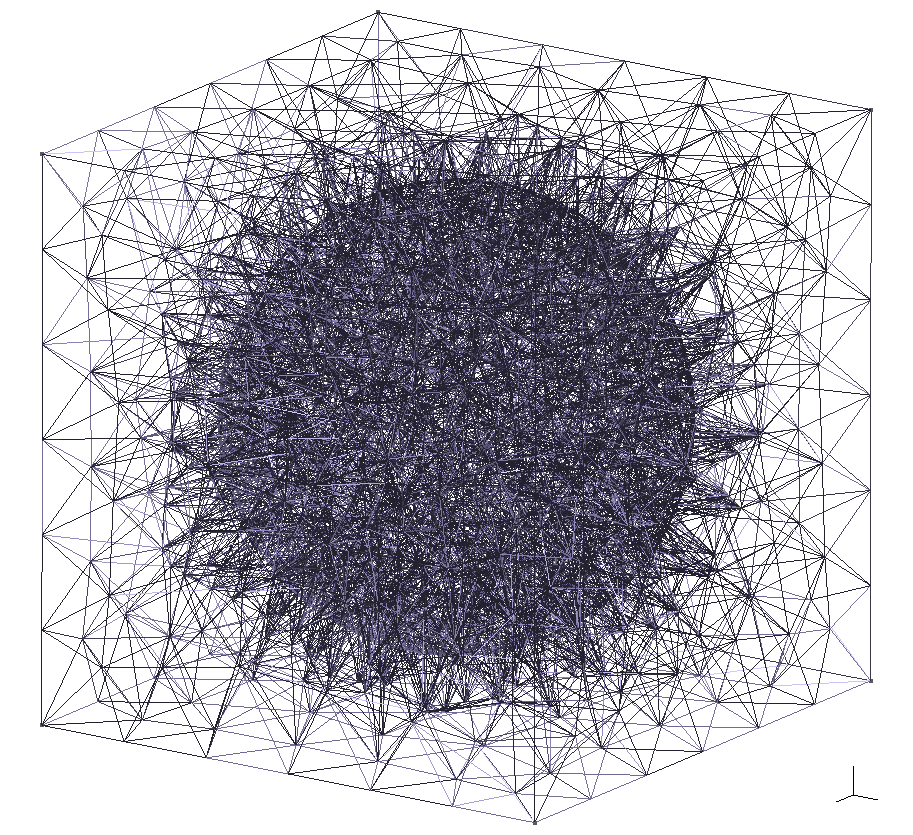}
    \end{center}
    \caption{\footnotesize
      Left: reference geometry with $4 \times 4 \times 4$ inclusions (only 1/8 of the geometry is used thanks to symmetry). Middle: macroscopic mesh. Right: mesh of the periodic cell.
    }
    \label{3D_HMM_Geo_mesh}
\end{figure}

\begin{figure}[H]
  \begin{center}
    \includegraphics[scale=0.14]{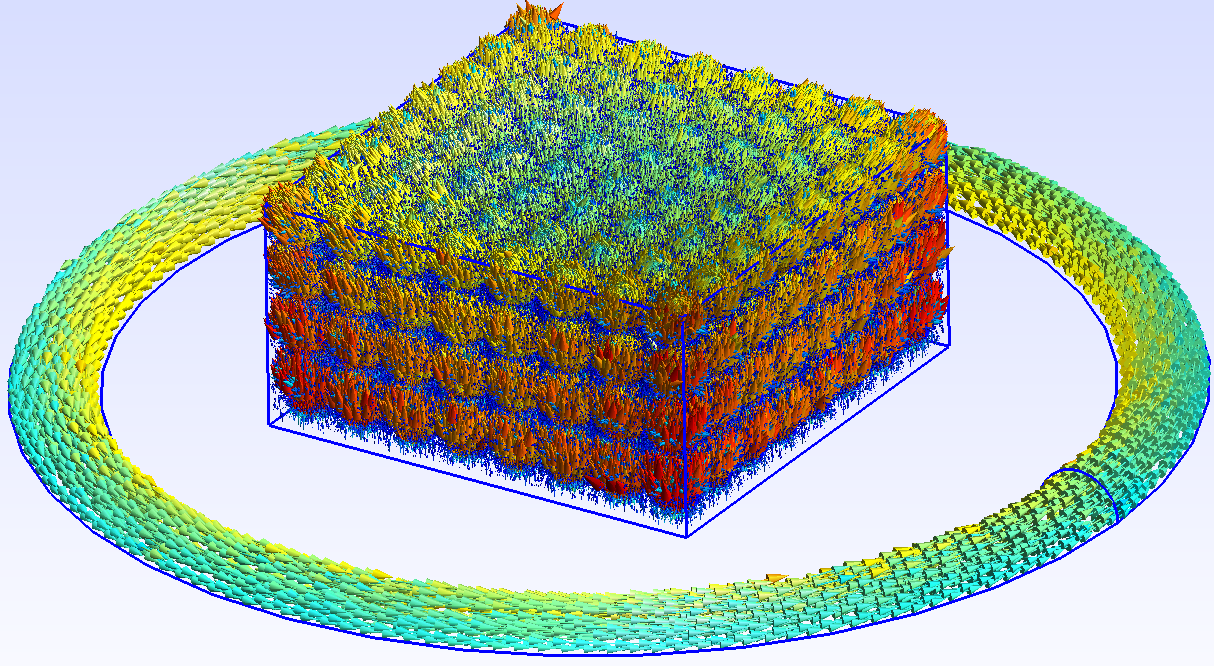}
  \includegraphics[scale=0.14]{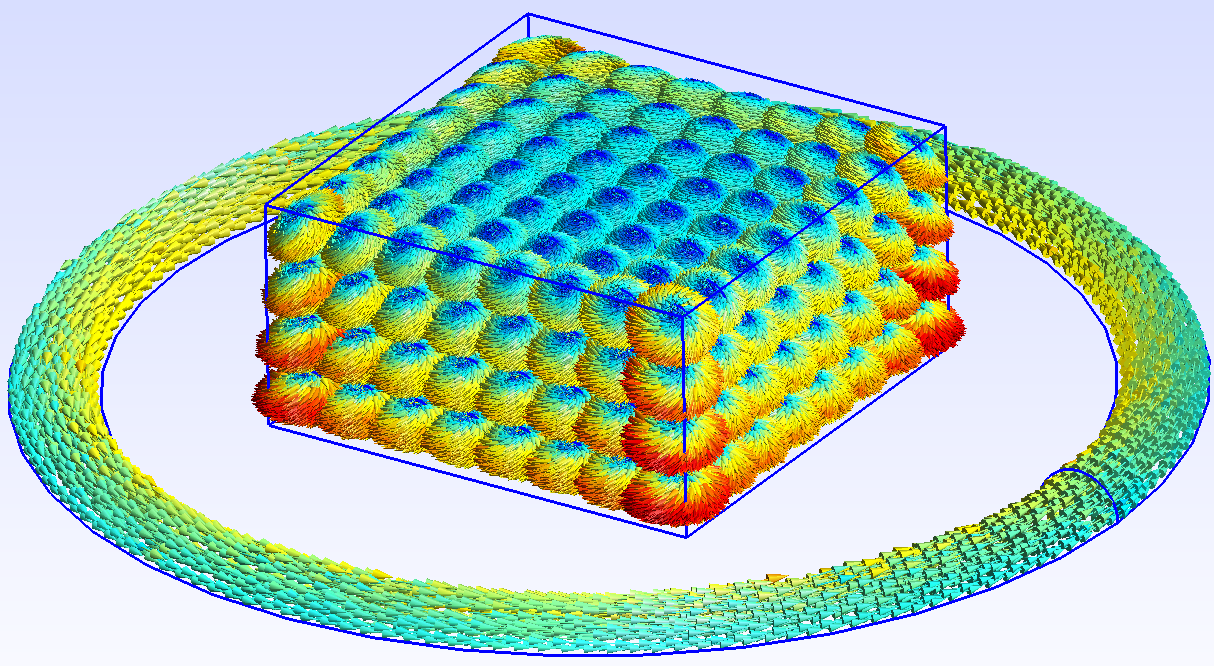}
  \end{center}
  \caption{\footnotesize
    Typical maps of the reference magnetic induction $\bb^{\varepsilon}$, the reference eddy currents $\bj^{\varepsilon}$ in the core, and the source current density $\bj_s$in the inductor at $\mu_r = 10^3$ and $f = 100$ kHz.
  }
  \label{HMM_3D_b_j_map}
\end{figure}

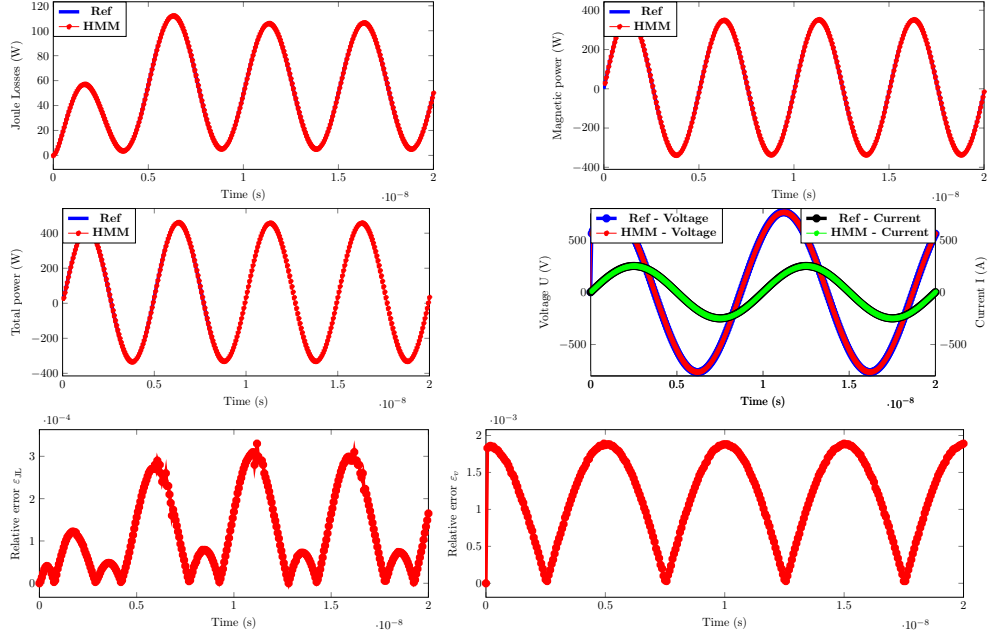
\begin{figure}[H]
  \begin{tikzpicture}[scale=0.45]
    \begin{axis}[xlabel={Time (s)}, ylabel={Joule Losses (W)}, xmin=0.0, xmax=2e-8, xtick={0, 0.5e-8, 1e-8, 1.5e-8, 2e-8}, width=0.98\columnwidth, height=0.5\columnwidth, legend style={at={(0,1)},anchor=north west}]
      \addplot [color=blue, mark=none, line width=1.0mm, mark size=2.0, mark options=solid] table {./data/3D_smc_320TS_lin_250A_Freq1e8_mur1e0_JouleLosses_Core_Ref.csv};
      \addlegendentry{\textbf{Ref}}
      \addplot [color=red, mark=*, mark size=2.0, mark options=dashdotted] table {./data/3D_smc_320TS_lin_250A_Freq1e8_mur1e0_Homog_JouleLosses.csv};
      \addlegendentry{\textbf{HMM}}
    \end{axis}
  \end{tikzpicture}
  \begin{tikzpicture}[scale=0.45]
    \begin{axis}[xlabel={Time (s)}, ylabel={Magnetic power (W)}, xmin=0.0, xmax=2e-8, xtick={0, 0.5e-8, 1e-8, 1.5e-8, 2e-8}, width=0.98\columnwidth, height=0.5\columnwidth, legend style={at={(0,1)},anchor=north west}]
      \addplot [color=blue, mark=none, line width=1.0mm, mark size=2.0, mark options=solid] table {./data/3D_smc_320TS_lin_250A_Freq1e8_mur1e0_MagPower_Core_Ref.csv};
      \addlegendentry{\textbf{Ref}}
      \addplot [color=red, mark=*, mark size=2.0, mark options=dashdotted] table {./data/3D_smc_320TS_lin_250A_Freq1e8_mur1e0_Homog_MagPower.csv};
      \addlegendentry{\textbf{HMM}}
    \end{axis}
  \end{tikzpicture}
  \begin{tikzpicture}[scale=0.45]
    \begin{axis}[xlabel={Time (s)}, ylabel={Total power (W)}, xmin=0.0, xmax=2e-8, xtick={0, 0.5e-8, 1e-8, 1.5e-8, 2e-8}, width=0.95\columnwidth, height=0.5\columnwidth, legend style={at={(0,1)},anchor=north west}]
      \addplot [color=blue, mark=none, line width=1.0mm, mark size=2.0, mark options=solid] table {./data/3D_smc_320TS_lin_250A_Freq1e8_mur1e0_TotPower_Core_Ref.csv};
      \addlegendentry{\textbf{Ref}}
      \addplot [color=red, mark=*, mark size=2.0, mark options=dashdotted] table {./data/3D_smc_320TS_lin_250A_Freq1e8_mur1e0_Homog_TotPower.csv};
      \addlegendentry{\textbf{HMM}}
    \end{axis}
  \end{tikzpicture}
  \begin{tikzpicture}[scale=0.45]
    \begin{axis}[xlabel={Time (s)}, ylabel={Voltage U (V)}, axis y line*=left, xmin=0.0, xmax=2e-8, xtick={0, 0.5e-8, 1e-8, 1.5e-8, 2e-8}, ymin=-800, ymax=800, width=0.90\columnwidth, height=0.5\columnwidth, legend style={at={(0,1)},anchor=north west}]
      \addplot [color=blue, mark=*, line width=1.0mm, mark size=2.0, mark options=solid] table {./data/3D_smc_320TS_lin_250A_Freq1e8_mur1e0_U.csv};
      \addlegendentry{\textbf{Ref - Voltage}}     
      \addplot [color=red, mark=*, mark size=2.0, mark options=dashdotted] table {./data/3D_smc_320TS_lin_250A_Freq1e8_mur1e0_U_Macro.csv};
      \addlegendentry{\textbf{HMM - Voltage}}
    \end{axis}
    \begin{axis}[xlabel={Time (s)}, ylabel={Current I (A)}, axis y line*=right, xmin=0.0, xmax=2e-8, xtick={0, 0.5e-8, 1e-8, 1.5e-8, 2e-8}, ymin=-800, ymax=800, width=0.9\columnwidth, height=0.5\columnwidth, legend style={at={(1,1)},anchor=north east}]
      \addplot [color=black, mark=*, line width=1.0mm, mark size=2.0, mark options=solid] table [x expr=\thisrowno{0}*1, y expr=\thisrowno{1}*-1, col sep=space] {./data/3D_smc_320TS_lin_250A_Freq1e8_mur1e0_I.csv};
      \addlegendentry{\textbf{Ref - Current}}     
      \addplot [color=green, mark=*, mark size=2.0, mark options=dashdotted] table [x expr=\thisrowno{0}*1, y expr=\thisrowno{1}*-1, col sep=space] {./data/3D_smc_320TS_lin_250A_Freq1e8_mur1e0_I_Macro.csv};
      \addlegendentry{\textbf{HMM - Current}}
    \end{axis}
  \end{tikzpicture}
  \begin{tikzpicture}[scale=0.45]
    \begin{axis}[xlabel={Time (s)}, ylabel={Relative error $\varepsilon_{\mathrm{JL}}$}, xmin=0.0, xmax=2e-8, xtick={0, 0.5e-8, 1e-8, 1.5e-8, 2e-8}, width=1.0\columnwidth, height=0.5\columnwidth, legend style={at={(0,1)},anchor=north west}]
      \pgfplotstableread{./data/3D_smc_320TS_lin_250A_Freq1e8_mur1e0_JouleLosses_Core_Ref.csv}{\datatable}
      \pgfplotstablecreatecol[copy column from table={./data/3D_smc_320TS_lin_250A_Freq1e8_mur1e0_Homog_JouleLosses.csv}{1}] {dataA} {\datatable}
      \addplot [color=red, mark=*, line width=1.0mm, mark size=2.0, mark options=solid] table [x expr=\thisrowno{0}*1e0, y expr=abs( (\thisrowno{1}/112 - \thisrow{dataA}/112))] {\datatable};
    \end{axis}
  \end{tikzpicture}
  \begin{tikzpicture}[scale=0.45]
    \begin{axis}[xlabel={Time (s)}, ylabel={Relative error $\varepsilon_{v}$}, xmin=0.0, xmax=2e-8, xtick={0, 0.5e-8, 1e-8, 1.5e-8, 2e-8}, width=1.20\columnwidth, height=0.5\columnwidth, legend style={at={(0,1)},anchor=north west}]
      \pgfplotstableread{./data/3D_smc_320TS_lin_250A_Freq1e8_mur1e0_U.csv}{\datatable}
      \pgfplotstablecreatecol[copy column from table={./data/3D_smc_320TS_lin_250A_Freq1e8_mur1e0_U_Macro.csv}{1}] {dataA} {\datatable}
      \addplot [color=red, mark=*, line width=1.0mm, mark size=2.0, mark options=solid] table [x expr=\thisrowno{0}*1e0, y expr=abs( (\thisrowno{1}/763 - \thisrow{dataA}/763))] {\datatable};
    \end{axis}
  \end{tikzpicture}
  \caption{\footnotesize 
    Global quantities. 
    \textbf{Top left}: Joule losses obtained solving the reference problem (blue curve) and the homogenized problem (red curve). 
    \textbf{Top right}: Magnetic power obtained solving the reference problem (blue curve) and the homogenized problem (red curve). 
    \textbf{Middle left}: Total electromagnetic power obtained solving the reference problem (blue curve) and the homogenized problem (red curve). 
    \textbf{Middle right}: Reference (blue curve) and homogenized (red curve) voltages, and reference (black curve) and homogenized (green curve) currents.
    \textbf{Bottom left}: Relative error $\varepsilon_{\mathrm{JL}}$ on Joule losses. 
    \textbf{Bottom right}: Relative error $\varepsilon_{v}$ on voltage. 
  }
  \label{fig:HMM_3D_lin_localquantities_1}
\end{figure}

\begin{figure}[H]
  \begin{tikzpicture}[scale=0.5]
    \begin{axis}[xlabel={Time (s)}, ylabel={Joule Losses (W)}, xmin=0.0, xmax=2.0e-5, xtick={0, 0.5e-5, 1.0e-5, 1.5e-5, 2.0e-5}, width=0.98\columnwidth, height=0.5\columnwidth, legend style={at={(0,1)},anchor=north west}]
      \addplot [color=blue, mark=none, line width=1.0mm, mark size=2.0, mark options=solid] table {./data/3D_smc_ref_fk_800A_f1e5_mur1e2_NbSteps160_JL_Core_Ref.csv};
      \addlegendentry{\textbf{Ref}}
      \addplot [color=red, mark=*, mark size=2.0, mark options=dashdotted] table {./data/3D_smc_hmm_fk_800A_f1e5_mur1e2_NbSteps160_Homog_JouleLosses.csv};
      \addlegendentry{\textbf{HMM}}
    \end{axis}
  \end{tikzpicture}
  \begin{tikzpicture}[scale=0.5]
    \begin{axis}[xlabel={Time (s)}, ylabel={Magnetic power (W)}, xmin=0.0, xmax=2.0e-5, xtick={0, 0.5e-5, 1.0e-5, 1.5e-5, 2.0e-5}, width=0.98\columnwidth, height=0.5\columnwidth, legend style={at={(0,1)},anchor=north west}]
      \addplot [color=blue, mark=none, line width=1.0mm, mark size=2.0, mark options=solid] table {./data/3D_smc_ref_fk_800A_f1e5_mur1e2_NbSteps160_MagPower_Core_Ref.csv};
      \addlegendentry{\textbf{Ref}}
      \addplot [color=red, mark=*, mark size=3.0, mark options=dashdotted] table {./data/3D_smc_hmm_fk_800A_f1e5_mur1e2_NbSteps160_Homog_MagPower.csv};
      \addlegendentry{\textbf{HMM}}
    \end{axis}
  \end{tikzpicture}
  \begin{tikzpicture}[scale=0.5]
    \begin{axis}[xlabel={Time (s)}, ylabel={Total power (W)}, xmin=0.0, xmax=2.0e-5, xtick={0, 0.5e-5, 1.0e-5, 1.5e-5, 2.0e-5}, width=0.95\columnwidth, height=0.5\columnwidth, legend style={at={(0,1)},anchor=north west}]
            \addplot [color=blue, mark=none, line width=1.0mm, mark size=2.0, mark options=solid] table {./data/3D_smc_ref_fk_800A_f1e5_mur1e2_NbSteps160_TotPower_Core_Ref.csv};
            \addlegendentry{\textbf{Ref}}
      \addplot [color=red, mark=*, mark size=3.0, mark options=dashdotted] table {./data/3D_smc_hmm_fk_800A_f1e5_mur1e2_NbSteps160_Homog_TotPower.csv};
      \addlegendentry{\textbf{HMM}}
    \end{axis}
  \end{tikzpicture}
  \begin{tikzpicture}[scale=0.45]
    \begin{axis}[xlabel={Time (s)}, ylabel={Voltage U (V)}, axis y line*=left, xmin=0.0, xmax=2e-5, xtick={0, 0.5e-5, 1e-5, 1.5e-5, 2e-5}, ymin=-2000, ymax=2000, width=0.93\columnwidth, height=0.5\columnwidth, legend style={at={(0,1)},anchor=north west}]
      \addplot [color=blue, mark=*, line width=1.0mm, mark size=2.0, mark options=solid] table {./data/3D_smc_ref_fk_800A_f1e5_mur1e2_NbSteps160_U.csv};
      \addlegendentry{\textbf{Ref - Voltage}}     
      \addplot [color=red, mark=*, mark size=2.0, mark options=dashdotted] table {./data/3D_smc_hmm_fk_800A_f1e5_mur1e2_NbSteps160_U_Macro.csv};
      \addlegendentry{\textbf{HMM - Voltage}}
    \end{axis}
    \begin{axis}[xlabel={Time (s)}, ylabel={Current I (A)}, axis y line*=right, xmin=0.0, xmax=2e-5, xtick={0, 0.5e-5, 1e-5, 1.5e-5, 2e-5}, ymin=-2000, ymax=2000, width=0.93\columnwidth, height=0.5\columnwidth, legend style={at={(1,1)},anchor=north east}]
      \addplot [color=black, mark=*, line width=1.0mm, mark size=2.0, mark options=solid] table [x expr=\thisrowno{0}*1, y expr=\thisrowno{1}*-1, col sep=space] {./data/3D_smc_ref_fk_800A_f1e5_mur1e2_NbSteps160_I.csv};
      \addlegendentry{\textbf{Ref - Current}}     
      \addplot [color=green, mark=*, mark size=2.0, mark options=dashdotted] table [x expr=\thisrowno{0}*1, y expr=\thisrowno{1}*-1, col sep=space] {./data/3D_smc_hmm_fk_800A_f1e5_mur1e2_NbSteps160_I_Macro.csv};
      \addlegendentry{\textbf{HMM - Current}}
    \end{axis}
  \end{tikzpicture}
  \begin{tikzpicture}[scale=0.5]
    \begin{axis}[xlabel={Time (s)}, ylabel={Relative error $\varepsilon_{\mathrm{JL}}$}, xmin=0.0, xmax=2e-5, xtick={0, 0.5e-5, 1e-5, 1.5e-5, 2e-5}, width=0.97\columnwidth, height=0.5\columnwidth, legend style={at={(0,1)},anchor=north west}]
      \pgfplotstableread{./data/3D_smc_ref_fk_800A_f1e5_mur1e2_NbSteps160_JL_Core_Ref.csv}{\datatable}
      \pgfplotstablecreatecol[copy column from table={./data/3D_smc_hmm_fk_800A_f1e5_mur1e2_NbSteps160_Homog_JouleLosses.csv}{1}] {dataA} {\datatable}
      \addplot [color=red, mark=*, line width=1.0mm, mark size=2.0, mark options=solid] table [x expr=\thisrowno{0}*1e0, y expr=abs( (\thisrowno{1}/0.0317 - \thisrow{dataA}/0.0317))] {\datatable};
    \end{axis}
  \end{tikzpicture}
  \begin{tikzpicture}[scale=0.5]
    \begin{axis}[xlabel={Time (s)}, ylabel={Relative error $\varepsilon_{v}$}, xmin=0.0, xmax=2e-5, xtick={0, 0.5e-5, 1e-5, 1.5e-5, 2e-5}, width=0.97\columnwidth, height=0.5\columnwidth, legend style={at={(0,1)},anchor=north west}]
      \pgfplotstableread{./data/3D_smc_ref_fk_800A_f1e5_mur1e2_NbSteps160_U.csv}{\datatable}
      \pgfplotstablecreatecol[copy column from table={./data/3D_smc_hmm_fk_800A_f1e5_mur1e2_NbSteps160_U_Macro.csv}{1}] {dataA} {\datatable}
      \addplot [color=red, mark=*, line width=1.0mm, mark size=2.0, mark options=solid] table [x expr=\thisrowno{0}*1e0, y expr=abs( (\thisrowno{1}/1630 - \thisrow{dataA}/1630))] {\datatable};
    \end{axis}
  \end{tikzpicture}
  \caption{\footnotesize
    Global quantities. 
    \textbf{Top left}: Joule losses obtained solving the reference problem (blue curve) and the homogenized problem (red curve). 
    \textbf{Top right}: Magnetic power obtained solving the reference problem (blue curve) and the homogenized problem (red curve). 
    \textbf{Middle left}: Total electromagnetic power obtained solving the reference problem (blue curve) and the homogenized problem (red curve). 
    \textbf{Middle right}: Reference (blue curve) and homogenized (red curve) voltages, and reference (black curve) and homogenized (green curve) currents.
    \textbf{Bottom left}: Relative error $\varepsilon_{\mathrm{JL}}$ on Joule losses. 
    \textbf{Bottom right}: Relative error $\varepsilon_{v}$ on voltage. 
  }
  \label{fig:HMM_3D_nonlin_localquantities_1}
\end{figure}
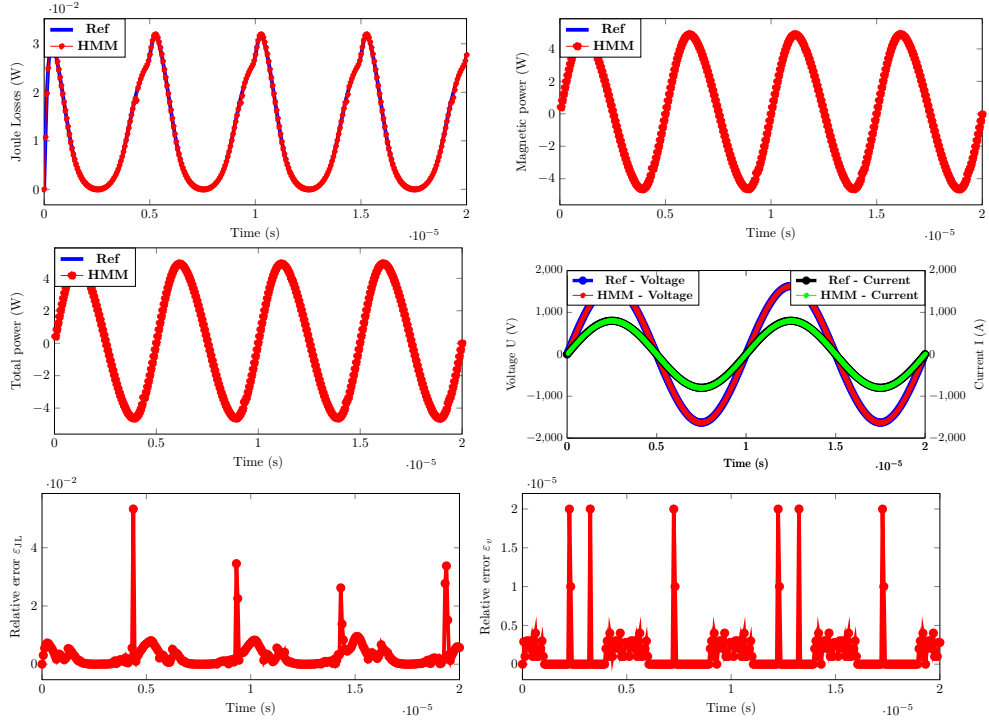

%
%
Reference results are obtained using a brute force approach with the FE problem solved on the fine mesh with 1,425,884 elements shown in Figure \ref{3D_HMM_Geo_mesh} (left). 
This mesh is obtained from a geometry of an idealized 3D SMC with $8 \times 8 \times 8$ magnetic and conducting balls, with the period $\Omega_m$ of width 100 $\mu$m and a conducting inclusion of radius 40 $\mu$m. Only 1/8$^{\mathrm{th}}$ of the geometry is considered thanks to symmetries.
HMM results are obtained using the macroscale mesh shown in Figure \ref{3D_HMM_Geo_mesh} (middle) with 64 macro-elements in the homogenized domain, and with 22,022 elements for the cell problem shown in Figure \ref{3D_HMM_Geo_mesh} (right).
Figure \ref{HMM_3D_b_j_map} shows typical maps of the reference magnetic induction $\bb^{\varepsilon}$ and eddy currents $\bj^{\varepsilon}$ in the core and the source current density $\bj_s$ in the inductor.
%
%
%
%
Similar to 2D problems, Figures~\ref{fig:HMM_3D_lin_localquantities_1} and \ref{fig:HMM_3D_nonlin_localquantities_1} present the time evolution of several global quantities (eddy-current losses, the magnetic power and the magnetic flux density $\bBB_M$ in the homogenized core, as well as the voltage), together with the corresponding the accuracy of these quantities as a function of the time. 

Figure \ref{fig:HMM_3D_lin_localquantities_1} shows results for a linear non-magnetic problem ($\mu_r = 1$) at $f = 100$ MHz. A sinusoidal source $I_{s}(t) = I_{s0} \sin{(2 \pi f t)}$ with amplitude $I_{s0} = 250$A is imposed and each period is discretized using 160 time steps. Figure \ref{fig:HMM_3D_nonlin_localquantities_1} shows results for a nonlinear problem with the Fr\"{o}hlich-Kennely nonlinear law with $\mu_{\mathrm{Fe}} = 10^2 \mu_0$ at $f = 100$ kHz. A sinusoidal source $I_{s}(t) = I_{s0} \sin{(2 \pi f t)}$ with amplitude $I_{s0} = 1000$A is imposed and each period is discretized using 160 time steps.

Figure \ref{fig:HMM_3D_lin_localquantities_1} shows a good agreement between the reference and HMM powers computed in the homogenized core : Joule losses (Top left), magnetic power (Top right) and the total power (Middle left)). The time-dependent voltage and current are depicted as well in Figure \ref{fig:HMM_3D_lin_localquantities_1} (Bottom-right). For this inductive system, the total eddy currents in the multiscale domain and in the massive inductor are comparable to the magnetic power and therefore one observes a phase shift between the current and the voltage (see the bottom right image in Figure 12). We also observe a good agreement between the reference and the homogenized voltages and currents from this image, with relative errors on Joule losses and on the voltage smaller than 0.035 \% for the Joule losses and smaller than 0.2 \% for the voltage. Figure \ref{fig:HMM_3D_nonlin_localquantities_1} also shows a good agreement for between the reference and the homogenized results.
%
%
%
%
\normalsize

The runtimes for the 2D and 3D simulations are reported in Table \ref{tab:runtime}. All simulations were performed on a computing cluster using a single CPU for the reference problems and 32 CPUs for the HMM simulations, except for the 3D HMM case, which was executed on 64 CPUs. A maximum wall-clock time of approximately two days was imposed by the cluster.

HMM simulations are not faster than the reference in these tests due to the small number of heterogeneous material periods in the reference problem. However, the computational complexity of HMM is independent from the number of periods, whereas the reference problem would become intractable on real world heterogeneous materials.
It should be noted that identical stopping criteria—a relative residual of $10^{-8}$ were employed for the reference and mesoscale problems. However, convergence at the macroscale is difficult to achieve, and even a small increase in the number of macroscale nonlinear iterations leads to a substantial rise in computational cost, since the mesoscale problem must be solved at each iteration.
\begin{table}[H]
  \centering
  \begin{tabular}{|c|c|c|c|c|}
    \hline
    Problem & 2D lin & 2D nonlin & 3D lin & 3D nonlin  \\
    \hline
    \hline
    REF     &  2h54  & 10h47     & 20h37  &  47h30      \\
    HMM     &  4h03  & 23h30     & 23h45  &  47h30      \\
    \hline
  \end{tabular}
  \caption{The simulation times for reference and the homogenized problems.}
  \label{tab:runtime}
\end{table}

%
%
%
\subsection{Three-dimensional SMCs: parallelization and performance of the HMM method}
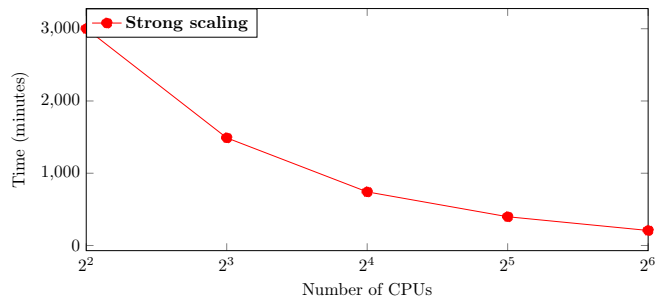
\begin{figure}[H]
    \centering
  \begin{tikzpicture}[scale=0.65]
    \begin{axis}[xlabel={Number of CPUs}, xmode=log, ylabel={Time (minutes)}, xmin=4.0, xmax=64.0, xtick={4, 8, 16, 32, 64}, log basis x = 2, width=1.0\columnwidth, height=0.5\columnwidth, legend style={at={(0,1)},anchor=north west} ]
      \addplot [color=red, mark=*, mark size=3.0, mark options=dashdotted] table {./data/3D_smc_StrongScaling_hmm_fk_800A_f1e5_mur1e2_NbSteps160.csv};
      \addlegendentry{\textbf{Strong scaling}}
    \end{axis}
  \end{tikzpicture}
  \caption{\footnotesize Strong scaling.}
  \label{fig:HMM_3D_nonlin_performance}
\end{figure}

%
%
Figure \ref{fig:HMM_3D_nonlin_performance} shows the results of strong scaling with a fixed macroscale mesh for the homogenized domain with increasing number of CPUs. As can be seen from this image, the cost of HMM decreases with the increasing number of CPUs.
\normalsize
%
%
\section{Conclusions and perspectives}
\label{sec:conclusions_perspectives}
%
%
This paper proposes an $\bh$-conforming multiscale formulation for magnetic problems with confined eddy currents. 
In the paper, two mesoscale problems are used to upscale homogenized quantities, a magnetoquasistatic problem, the solution of which is used for upscaling the homogenized magnetic flux density $\bBB_M$ and a magnetostatic problem, the solution of which is used for upscaling the macroscale incremental reluctivity $(\partial \bBB_M/\partial \bHH_M)$. Additionally, it is shown that the magnetoquasistatic problem can be replaced by a magnetostatic problem provided that the eddy currents are negligeable compared to the magnetic power in the cell.

Newton--Raphson algorithms are also used to speed up the convergence of the nonlinear problems at both scales with a relaxation factor $\omega \in [0.05, 1]$. Approximately 10-12 non linear iterations were necessary for the convergence of some macroscale time steps. For convergence of the macroscale nonlinear problem, the convergence criteria based on the relative increment of the solution and on physical powers (such as eddy current losses, magnetic power, or total power) demonstrated robustness compared to convergence criteria based on the residual $|\br| = |\bb - \bA \bx|$ from the macroscale linear system. Solutions obtained through the multiscale method demonstrated accuracy, with global quantities such as eddy current losses, magnetic power, and currents/voltages converging to values consistent with those obtained using the brute force approach. 

The convergence of the solutions and global quantities depends on various parameters, including the macroscale mesh resolution. The macroscale problems exhibit effective strong scaling with an increasing number of CPUs, resulting in reduced computational costs. However, weak scaling is not fully guaranteed, as the computational cost increases noticeably from 2 to 16 CPUs and stabilizes thereafter from 16 to 32 CPUs.

These findings represent a significant advancement towards modeling real-life 3D composites characterized by nonlinear magnetic laws and the presence of significant confined eddy currents within the cell. Areas for potential improvement include further investigation of the robustness of convergence of the methodology, reducing the computational cost by the evaluation of the homogenized incremental reluctivity using sensitivity analysis, and the use of model order reduction techniques to minimize the number of mesoscale problems solved at the macro level and simplify the reconstruction of mesoscale solutions. 
In the present work, the cell geometry was assumed to be periodic; extending the approach to handle multiscale materials with nonperiodic mesoscale geometries constitutes another important direction for future research.
%
%
%
%
%
%
\bibliographystyle{siamplain}
\bibliography{./3D_HMM_H_Conforming.bib}
\end{document}